
\magnification =\magstep1
\baselineskip =13pt
\overfullrule =0pt

\def\ldb{[\![}
\def\rdb{]\!]}
\def\LDB{\big[\!\!\big[}
\def\RDB{\big]\!\!\big]}

\def\bb{\bf}
\def\Bigldb{\left[\vbox  to 21pt{}\right.\kern-6pt\left[\vbox to
21pt{}\right.}
\def\Bigrdb{\left.\vbox to 21pt{}\right]\kern-6pt\left.\vbox to
21pt{}\right]}
\def\<{{\langle}}
\def\>{{\rangle}}
\def\[{{[\! [}}
\def\]{{]\!]}}

\centerline {\bf NONCOMMUTATIVE GEOMETRY BASED ON COMMUTATOR EXPANSIONS}

\vskip .7cm

\centerline {\bf 
M. Kapranov}

\vskip 1cm

\centerline {\bf Contents}

\vskip .4cm

{1. NC-complete algebras.

2. NC-schemes.

3. The NC-affine space and Feynman-Maslov operator calculus.

4. Detailed study of algebraic NC-manifolds.

5. Examples of NC-manifolds. }

\vskip .4cm

\noindent The term ``noncommutative geometry" has come to signify a vast framework
of ideas directed towards generalization, to noncommutative rings, of the fundamental duality between spaces and (commutative) rings of functions
on them. This duality is at the basis of the modern approach to geometry.
The main motivation to look for noncommutative generalizations is
provided by quantum mechanics where commuting functions are replaced by
non-commuting operators.

The study of noncommutative rings by methods and concepts inspired
by differential geometry such as (analogs) of differential forms,
connections, integrations,  Chern classes etc. was pursued very
actively and many important results have been obtained, including
the discovery of cyclic homology [Co]. 

From the point of view of algebraic geometry, however,
commutative rings by themselves correspond to only a particular class of
``spaces", namely affine schemes, and the thrust of the theory is that
they can be glued to form more global objects.  The construction
of such objects in the noncommutative case
(and even of any actual geometric objects
corresponding to noncommutative rings)  has proved to be
quite difficult despite many interesting developments [AZ] [Ros] [VV].

The aim of the present paper is to develop an approach to noncommutative
algebraic geometry ``in the perturbative regime" around ordinary
commutative geometry. Let $R$ be a noncommutative algebra (over {\bf C})
and $R_{ab}= R/[R,R]$ be its commutativization. Then we know the
geometric object $X_{ab}={\rm Spec}(R_{ab})$. The naive aim of
noncommutative algebraic geometry would be to associate to the surjection
$R\to R_{ab}$ an embedding of $X_{ab}$ into some ``noncommutative space" 
$X= {\rm Spec}(R)$. The essense of our perturbative approach is not to
worry about the whole $X$ but concentrate on the formal neighborhood
of $X_{ab}$ in $X$. This neighborhood can be described, as in the usual
algebro-geometric
theory of formal schemes, by equipping $X_{ab}$, a known object,
by an appropriate sheaf of noncommutative rings ${\cal O}^{NC}$. 
On the algebraic level, this means that we complete $R$ with respect to
the topology, in which iterated commutators, like
$[a_1, [a_2, ..., [a_{n-1}, a_n]...]$, or $[a_1, b_1] \cdot ... \cdot [a_n, b_n]$,
are small. Elements of the completion can be thought of as formal
commutator series.  

Further, ringed spaces of the form $(X_{ab}, {\cal O}^{NC})$ 
can be glued
together to form more general objects which we call NC-schemes.
Any NC-scheme $X$ consists of an ordinary scheme $X_{ab}$ and an
appropriate sheaf  of rings ${\cal O}^{NC}$ on it.

We are especially interested in NC-manifolds, which are NC-schemes 
$(M, {\cal O}^{NC})$ where $M$ is a smooth algebraic
variety of some dimension  $n$, and the completion
of ${\cal O}^{NC}$ at any $x\in M$ is isomorphic to
${\bf C}\<\!\< x_1, ..., x_n\>\!\>$, the completion of the
free associative algebra. This concept globalizes, in a sense, the theory of 
 smooth (quasi-free) noncommutative algebras
studied by W. Schelter, J. Cuntz and D. Quillen [Sche][Q].

 Given an ordinary manifold $M$, an NC-manifold
(thickening)  $X$ with $X_{ab}=M$
 can be thought of as a sophisticated
differential-geometric structure on $M$, leading to additional
characteristic classes, e.g., the NC-Atiyah class in $H^1(M, \Omega^2_M\otimes
T_M)$, see (4.5).
The idea that one can and should seriously develop noncommutative geometry
based on free associative algebras (rather than on algebras
with relations resembling or relaxing commutativity) was put forward
by M. Kontsevich [Ko] in the formal  case and I.M. Gelfand and
V.S. Retakh [GR] in the affine case. In fact, M. Kontsevich communicated
to the author that he was aware of the possibility of a theory such 
as developed here. 

\vskip .3cm

The paper is organized as follows. In Section 1 we study NC-nilpotent
algebras, i.e., those in which higher iterated commutators vanish. 
In particular, we introduce the concept of a $d$-smooth algebra
and show (Theorem 1.6.1) that any smooth finitely generated commutative algebra
admits a unique (up to a non-canonical isomorphism) $d$-smooth extension.
The relation between  smooth commutative algebras and  $d$-smooth algebras 
is somewhat similar to the relation between free abelian groups and
free $d$-stage nilpotent groups. The class of NC-smooth algebras
obtained by passing to the limit $d\to\infty$, is analogous to
the class of free pro-nilpotent groups. 

In Section 2 we describe the geometric objects (affine NC-schemes)
corresponding to NC-nilpotent
algebras as well as well as more global objects obtained by gluing these
ones. In particular, we single out the class of algebraic NC-manifolds
(whose algebras of functions are NC-smooth) and prove (Theorem 2.3.5)
that such manifolds can be identified with appropriate functors on
the category of NC-nilpotent algebras. 

Section 3 is devoted to a detailed study of the simplest NC-manifold,
the noncommutative affine space. It is represented by a certain
sheaf ${\cal O}^{NC}$ of noncommutative rings on the usual affine
space. We give a completely  explicit description
(Theorem 3.5.3) of this sheaf by
using the ideas from the Feynman-Maslov ``calculus of ordered operators". 
The appearance of iterated derivatives in the explicit formulas gives
some insight into the differential-geometric meaning of noncommutativization. 

In Section 4 we study the problem of constructing a NC-thickeninng
of a given  ordinary manifold, by homological
means. We exhibit a series of cohomological obstructions whose
vanishing is necessary and sufficient for the existence of a thickening. 
These obstructions have a nice interpretation in terms of 
 $A_\infty$-structures of J. Stasheff [Sta], the very first obstruction
being a certain associator. 

Finally, in Section 5 we show that several familiar algebraic varieties
possess natural NC-thickenings. These include all the classical flag varieties
and all the smooth moduli spaces of vector bundles.
In particular,
for the case of Grassmannians what we get can be seen as a 
completion, along the commutative points, of the 
noncommutative Grassmannians  of Gelfand-Retakh [GR]
which are described as certain explicit functors on the category of
skew fields defined in terms of analogs of affine charts.  

\vskip .2cm

This research was partially supported by an NSF grant.

\vfill\eject

\centerline {\bf \S1. NC-complete algebras.}

\vskip 1cm

\noindent
{\bf (1.1) The NC-filtration.} In this paper the base
field is the field $\bb C$ of complex numbers. In particular, all algebras
are $\bb C$- algebras.

Let $L$ be a ${\rm Lie}$ algebra.  Its lower central series is the sequence
of subalgebras

$$L_m = [L, [L, \ldots, [L,L]\dots]\qquad \hbox{($m$ times),}\quad m \ge 1$$

\noindent
Thus $L_1=L$ and $L_m$ is spanned by the expressions $[x_1,
[x_2,\ldots[x_{m-1},x_m]\ldots]$ containing $m-1$ instances of Lie
brackets.

Let now $R$ be an associative algebra and $R^{\rm Lie} = (R, [a,b] = ab -
ba)$ be $R$ regarded as a ${\rm Lie}$ algebra.

\proclaim (1.1.1) Definition.  The NC-filtration of $R$ is the
decreasing filtration $\{F^dR\}_{d\ge 0}$ where $F^dR$ is the two-sided
ideal
$$F^dR = \sum_m \sum_{i_1+\cdots+ i_m -m = d} R \cdot R^{\rm Lie}_{i_1}
\cdot R \cdot{\kern 1.5pt}{\cdots}{\kern 1.5pt}\cdot R\cdot R^{\rm
Lie}_{i_m}\cdot R$$

\noindent
Thus $F^0R=R$ and $F^dR$ is generated by  expressions containing $d$
instances of commutator brackets. This filtration was considered by Helton
and Howe [HH] under the name ``commutator filtration.'' Its main property is
as follows.

\proclaim (1.1.2) Proposition. $F$ is a decreasing algebra
filtration (i.e., $(F^{d_1}R)(F^{d_2}R)\subset F^{d_1+d_2}R)$ such that the
associated graded algebra $\hbox{gr}^{\bullet}_F R = \bigoplus F^dR/F^{d+1}R$
is commutative. In other words, $[F^{d_1}R, F^{d_2}R] \subset
F^{d_1+d_2+1}R$.

The number

$$\hbox{ord}_{NC}(f) =\min \{d \colon\  f\in F^dR\}$$

\noindent
will be called the NC-order of $f\in R$. Note that $\hbox{gr}^0_F(R) =
R/[R,R]$ is the commutativization of $R$. This algebra will be denoted by
$R_{ab}$. Every $\hbox{gr}^i_F(R)$ is an $R_{ab}$-module. We have a
canonical surjective homomorphism $R\rightarrow R_{ab}$ whose value on $f\in
R$ will be denoted $f_{ab}$.

We define the NC-topology on $R$ to be the topology in which the $F^dR$ form
a basis of neighborhoods of $0$.
Recall that a Poisson algebra is a commutative algebra $P$ equipped with an
anticommutative binary operation $\{f,g\}$, called the Poisson bracket,
which satisfies the Jacobi identy and is a derivation with respect to each
argument. A Poisson algebra $P$ will be called graded, if $\displaystyle
P=\mathop{{}\oplus{}}_{d\in {\bb Z}}P^d$ is ${\bb Z}$-graded as a vector
space, and
$$P^{d_1}\cdot P^{d_2} \subset P^{d_1
+d_2}, \quad \{P^{d_1}, P^{d_2}\}\subset P^{d_1+d_2+1}.\leqno (1.1.3)$$

\noindent Proposition  1.1.2 implies, in a standard way, the following.

\proclaim (1.1.4) Proposition. The algebra
$\hbox{gr}^{\bullet}_F(R)$ has a natural structure of a graded Poisson
algebra.

\proclaim (1.1.5) Definition. (a) An associative algebra $R$
is called NC-nilpotent of degree $d$ (resp. NC-nilpotent), if $F^{d+1}R=0$
(resp. $F^iR=0$ for $i\gg 0$).
\hfill\break
(b) For an associative algebra $R$ its NC-completion is the algebra
$$R_{\ldb ab \rdb} = \lim_{\leftarrow} R/F^dR$$
(c) $R$ is called NC-complete, if the natural morphism $R \rightarrow
R_{\ldb ab \rdb}$ is an isomorphism.

\proclaim (1.1.6) Proposition.   The following are equivalent:\hfill\break
(i) $R$ is NC-complete.\hfill\break
(ii) $R$ is Hausdorff and complete with respect to the NC-topology.
\hfill\break
{(iii)} $R$ is an inverse limit of NC-nilpotent algebras.

\vskip.2cm

\noindent {\bf (1.2) Central extensions.} Let $R$ be an associative
albebra, $M$ be an $R$-bimodule. We say that $M$ is central, if $am = ma$
for any $a\in R, m\in M$.

\proclaim (1.2.1) Proposition. A central $R$-bimodule is the
same as a module over $R_{ab}$.

\proclaim (1.2.2) Definition.  Let $R$ be an associative
algebra. An Abelian extension of $R$ is an exact sequence of algebras
$$0\rightarrow I\rightarrow R^{\prime}\rightarrow R\to 0$$
where the ideal $I$ satisfies $I^2 = 0$. A central extension is an Abelian
extension such that $I$ lies in the center of $R^{\prime}$.

Thus, in an Abelian extension, the only nontrivial structure on $I$ is that
of an $R$-bimodule. For a central extension, $I$ is a central $R$-bimodule.

\proclaim (1.2.3) Proposition. (a) NC-nilpotent algebras are
precisely algebras which can be obtained as iterated central extensions of
commutative algebras. \hfill\break
(b) Any surjection $R_1 \rightarrow R_2$ of NC-nilpotent algebras whose
kernel is a nilpotent ideal, can be decomposed into a sequence of central
extensions.

Given an algebra $R$ and an $R$-bimodule $M$, we will denote by $R \oplus M$
the trivial Abelian extension of $R$ by $M$, which is the direct sum with
the product given by
$$(a_1, m_1)(a_2, m_2) = (a_1a_2, m_2a_2 +
a_1 m_2).\leqno (1.2.4)$$

\proclaim (1.2.5) Proposition. Let
$$0\rightarrow I\rightarrow R^{\prime} \mathop{{}\rightarrow{}}^{p} R
\rightarrow 0$$
be a central extension. Then:\hfill\break
(a) The group of automorphisms of the extension identical on $I, R$ is
Abelian and is identified with $\hbox{Der}(R_{ab}, I)$, the module of
$I$-valued derivations.\hfill\break
(b) The fiber product $R^{\prime} \mathop{{}\times{}}_R R^{\prime}$ is
identified with
$$R^{\prime} \mathop{{}\times{}}_{R_{ab}} (R_{ab} \oplus I).$$

\noindent {\sl Proof;} (a) Given any Abelian (not necessarily central) extension,
the group of its automorphisms as above is identified with
$$\hbox{Der}(R,I) = \{D \colon\  R\rightarrow I\mid D(ab) = a\, D(b) +
D(a)\, b\}.$$
Indeed, for such an automorphism $g$ the map $g-1 \colon\  R'\rightarrow R'$
takes values in $I$ and vanishes on $I$, so descends to $D\colon\
R\rightarrow I$, which lies in $\hbox{Der}(R,I)$.

Now, since $I$ is central, one verifies by the Leibniz rule that $D(abcd) =
D(acbd)$, so $D$ descends to a derivation  $R_{ab}\rightarrow I$. This shows
that $\hbox{Der}(R,I) = \hbox{Der}(R_{ab}, I)$.

(b) Define
$${\varphi}\colon\  R^{\prime} \mathop{{}\times{}}_{R} R^{\prime}
\to R'
\mathop{{}\times{}}_{R_{ab}} (R_{ab} \oplus I), 
\quad (x,y)\mapsto (x,x_{ab} +
y -x).$$
It is clear that $\varphi$ is bijective. The fact that $\varphi$ is a
homomorphism follows because $I^2 = 0$ and the $R$-action on $I$ factors
through $R_{ab}$.

\vskip .1cm

\proclaim  (1.2.6) Proposition. Let $p\colon\ R^{\prime}
\rightarrow R$ be a central extension as in  (1.2.5) and $g\colon\
R^{\prime} \rightarrow R^{\prime}$ an algebra endomorphism such that $pg=p$.
If $p_*\colon\ R^{\prime}_{ab} \rightarrow R_{ab}$ is an isomorphism, then
$g\big\vert_I = Id$. In particular, $g$ is an isomorphism.

\noindent {\sl Proof:} Consider the homomorphism
$$\Psi\colon\ R^{\prime} \mathop{{}\rightarrow{}}^{(Id,g)} R^{\prime}
\mathop{{}\times{}}_R R^{\prime} \mathop{{}\simeq{}}^{\varphi} R^{\prime}
\mathop{{}\times{}}_{R_{ab}} (R_{ab} \oplus I) \rightarrow (R_{ab} \oplus
I).$$
The algebra $R_{ab} \oplus I$ is commutative, so this homomorphism descends
to a homomorphism $\psi: R_{ab}^{\prime}\rightarrow
(R_{ab} \oplus I)$ whose projection to $R_{ab}$ is $p_*$. Since $p_*$ is an
isomorphism, we can view $\psi$ as a homomorphism $R_{ab}\rightarrow
R_{ab}\oplus I$ whose composition with the projection to $R_{ab}$ is
the identity. The set of such homomorphisms is identified with
$\hbox{Der}(R_{ab}, I)$. From the nature of the identification ${\varphi}$
it follows that $g(x) = x + D(p(x))$ for some $D\in \hbox{Der}(R_{ab}, I)$,
so $g$ is an automorphism of the form described in (1.2.5)(a).

\vskip .1cm

\proclaim (1.2.7) Corollary. Let $R^{\prime}_i
\buildrel p_i\over\rightarrow R, i=1,2$, be two central extensions of $R$
with $p_{i*}\colon\ R^{\prime}_{i,ab} \rightarrow R_{ab}$ isomorphisms.
Suppose there are homomorphisms $f\colon\ R^{\prime}_1 \rightarrow
R^{\prime}_2$, $g\colon\ R^{\prime}_2 \rightarrow R^{\prime}_1$ compatible
with the $p_i$. Then $R^{\prime}_1$ is isomorphic to $R^{\prime}_2$.
Moreover, the kernels of the $p_i$ are identified in a canonical way (i.e.,
the identification is independent of the choice of $f,g$, provided they
exist).

\vskip .2cm

\noindent {\bf (1.3) Hochschild homology and universal central
extensions.}
Let $R$ be an associative algebra, $I$ be an $R$-bimodule. We denote by
$H_{\bullet}(R,I), H^{\bullet}(R,I)$ the Hochschild homolohy and cohomology
of $R$ with coefficients in $I$, see $[L]$. The corresponding chain and
cochain complexes have the form 
$$ C_{\bullet}(R,I) =
\{I\mathop{{}\otimes{}}_{\bb C} R^{\otimes n}, n\geq 0\}, C^{\bullet}(R,I) =
\{\hbox{Hom}_{\bb C}(R^{\otimes n},I), n\geq 0\}.$$  The following is well
known.

\proclaim (1.3.1) Proposition. $H^2(R,I)$ is identified with
the set of Abelian extensions (1.2.2) modulo isomorphisms identical on $R$
and $I$.

\noindent We now concentrate on the case when $I$ is central.

\proclaim (1.3.2) Proposition. If $I$ is a central
$R$-bimodule, then $C_{\bullet}(R,I), C^{\bullet}(R,I)$ have natural
structures of complexes of $R_{ab}$-modules. In particular, $H_m(R,I),
H^m(R,I)$ are $R_{ab}$-modules.

\noindent {\sl Proof:} Consider, for example, the chain complex. We
introduce the $R_{ab}$-module structure in $C_m(R,I) =
I\mathop{{}\otimes{}}\limits_{\bb C} R^{\otimes m}$ by viewing $R^{\otimes
n}$ as the vector space of multiplicities for the $R_{ab}$-module $I$, i.e.,
$$a\otimes (i\otimes b_1 \otimes\cdots \otimes b_m)=ai\otimes
b_1\otimes\cdots \otimes b_m, a\in R_{ab}.$$
The boundary in $C_{\bullet}(R,I)$ has the form
$$\eqalign{
\partial (i\otimes b, \otimes \cdots \otimes b_m) ={}& ib_1 \otimes
b_2\otimes \cdots \otimes b_m + \sum^{m-1}_{\nu =1} (-1)^\nu i\otimes \cdots
\otimes b_\nu b_{\nu +1} \otimes\cdots\otimes b_m\cr
&{}+ (-1)^m b_mi\otimes b_1\otimes \cdots \otimes b_{m-1}\cr}.$$
Because $I$ is central, $\partial$ commutes with the action of $a$.
Proposition is proved.

\vskip.1cm

The ``universal'' example of a central $R$-bimodule is given by $R_{ab}$.
The following fact is well known if $R=R_{ab}$ is commutative ($[L]$, Prop.
1.1.10).

\proclaim  (1.3.3) Proposition. Let $R$ be any associative
algebra. Then $H_1(R,R_{ab}) \simeq \Omega^1_{R_{ab}}$ is the module of
K\"ahler differentials of $R_{ab}$.

\noindent {\sl Proof:} This is a modification of the proof in
{\it loc.cit.} The relevant part of the Hochschild complex is
$$\displaylines{
R_{ab} \otimes R\otimes R\mathop{{}\rightarrow{}}^{\partial_2} R_{ab}
\otimes R\mathop{{}\rightarrow{}}^{\partial_1} R_{ab},\cr
\partial_1(f\otimes g) = f\cdot G_{ab} - G_{ab}\cdot f=0,\cr
\partial_2 (f\otimes G\otimes H) = f\otimes GH - fG_{ab} \otimes H -
H_{ab}f\otimes G.\cr}$$
Thus $H_1(R,R_{ab}) = \hbox{Coker}(\partial_2)$. We now define the maps
$$\displaylines{
{\varphi}\colon\ H_1 (R,R_{ab}) \rightarrow \Omega^1_{R_{ab}}, 
\quad f\otimes
G\mapsto f\cdot d(G_{ab}),\cr
\psi\colon\ \Omega^1_{R_{ab}}\rightarrow H_1(R,R_{ab}),
\quad  f\, dg \rightarrow
f\otimes G\quad \hbox{mod}\quad \hbox{Im}(\partial_2).\cr}$$
Here $G\in R$ is any element such that $G_{ab} = g$. The well-definedness of
${\varphi}$ is clear. To see that $\psi$ is well defined, we need to prove
that for any $A,G_1,G_2,B\in R, F\in R_{ab}$ we have
$$f\otimes AG_1G_2B\equiv f\otimes AG_2G_1B\quad \hbox{(mod
Im($\partial_2$))}.$$
To see this, denote $a=A_{ab}$, $g_i = G_{i,ab}$ etc. and find:
$$f\otimes AG_1G_2B\equiv fag_1g_2\otimes B + bfag_1\otimes G_2 + g_1g_2bf
\otimes A+ g_2bfa\otimes G_1$$.
Since $g_2,g_2$ commute in $R_{ab}$, the right hand side for $f\otimes
AG_2G_1B$ will be the same. Having established the existence of $\varphi$
and $\psi$, it is immediate that  they are mutually inverse.

\vskip .1cm

We now note the following universal coefficient formula.

\proclaim (1.3.4) Proposition. Let $R$ be any associative
algebra, $I$ a central $R$-bimodule. Then there is a spectral sequence
$$E^{ij}_2 = \hbox{Ext}^j_{R_{ab}} (H_i(R,R_{ab}), I) \Rightarrow
H^{i+j}(R,I).$$

\noindent {\sl Proof:} This follows from the identification of
complexes
$$C^{\bullet}(R,I)\cong \hbox{Hom}_{R_{ab}}(C_{\bullet}(R,R_{ab}),I)$$
and from the fact that $C_{\bullet}(R,R_{ab})$ consists of free
$R_{ab}$-modules.

\proclaim (1.3.5) Corollary. Let $R$ be an associative
algebra such that $R_{ab}$ is smooth. Then
$$H^2(R,I)=\hbox{Hom}_{R_{ab}}(H_2(R,R_{ab}),I)$$
for a any central $R$-bimodule $I$.
In particular, in this situation we have the tautological class
$$\tau_R \in H^2(R,H_2(R,R_{ab}))$$
corresponding to the identity map of $H_2(R,R_{ab})$.

\noindent {\sl Proof:} Since $R_{ab}$ is smooth,
$H_1(R,R_{ab})=\Omega^1_{R_{ab}}$ is projective, so
$\hbox{Ext}^j_{R_{ab}}(H_1(R,R_{ab}), I) =0, j>0$.
Also $\hbox{Ext}^j_{R_{ab}}(H_0(R, R_{ab}),I) =
\hbox{Ext}^j_{R_{ab}}(R_{ab},I)=0, \quad j> 0.$ Thus the only nontrivial
term $E^{ij}_2$ with $i+j=2$ is $\hbox{Hom}_{R_{ab}}(H_2(R,R_{ab}),I)$. There
are no differentials hitting this term, and the only possible differentials
originating from it, are  $d_2$, with values in
$\hbox{Ext}^2_{R_{ab}}(H,(R,R_{ab}),I)=0$, and $d_3$, with values in a
subquotient of $\hbox{Ext}^3_{R_{ab}}(H_0(R,R_{ab})),I)=0$. This proves our
assertion.

\proclaim (1.3.6) Definition. Let $R$ be an associative
algebra such that $R_{ab}$ is smooth. The universal central extension of $R$
is the extension
$$0\rightarrow  H_2(R,R_{ab})\rightarrow R^{\tau}\rightarrow R\rightarrow
0$$
corresponding to the  tautological Hochschild class $\tau_R$.

Since we have only the class $\tau_R$ but  not, in general, a distinguished
cocycle representing it, we need to make precise in which sense we can speak
about ``the'' universal central extension.

\proclaim
 (1.3.7) Proposition. (a) The universal central
extension $R^{\tau}$ is defined uniquely up to a (non-canonical) isomorphism
identical on $R$, $H_2(R,R_{ab})$.\hfill\break
(b) Let $U(R)$ be the category (groupoid) formed by universal central
extensions and their isomorphisms as in (a). The $U(R)$ is a gerbe with band
$\hbox{Der}(R_{ab}, H_2(R,R_{ab}))$.

Part (b) means that for any two choices $R^{\tau}_1$ and $R^{\tau}_2$ of the
universal central extension,  the set of isomorphisms $R^{\tau}_1\rightarrow
R^{\tau}_2$ identical on $R$, $H_2(R,R_{ab})$ is  a principal homogeneous
space over $\hbox{Der}(R_{ab},H_2(R,R_{ab}))$, and this structure of
principal homogeneous space is compatible with the composition of
isomorphisms, see [Bry]. This property follows from Proposition 1.2.5(a).

\proclaim (1.3.8) Proposition. Let $R$ be as above and
$$0\rightarrow J\rightarrow S\rightarrow R\rightarrow 0$$
be any central extension. Then there is a morphism of extensions $\gamma\colon\
R^{\tau} \rightarrow S$, $\gamma (H_2(R,R_{ab})) \subset J$, identical on $R$.
This morphism
 has the property that for any morphism of extensions $\chi\colon\
S  \rightarrow R^{\tau}$, $\chi(J) \subset H_2(R, R_{ab})$, the
map $\chi\gamma$ is identical on $H_2(R, R_{ab})$.

\noindent {\sl Proof:} Let $H=H_2 (R, R_{ab})$ and $c \in H_2
$ be the class corresponding to $S$. By (1.3.5) we can identify $c$ with a
morphism $\tilde{c}\colon\ H
\rightarrow J$. Such a morphism defines a morphism $\gamma$ of extensions, as
claimed. Let now $\chi$ be given. Then $\chi \big\vert_J\colon\ J  \rightarrow
H$ is a morphism of left $R$-bimodules. The fact that $\chi$ is a morphism of
extensions means that the image of $c$ under $(\chi\big\vert_J)_*\colon\ H^2
(R,J)  \rightarrow H^2(R,H)$ is $\tau_R$. But this is  the same as saying
that the composition $\chi\gamma\colon\ H  \rightarrow H$ is the identity.
\vskip .2cm

\noindent
{\bf\ (1.3.9) Example.} Let $R$ be a smooth commutative
algebra. Then $R_{ab} = R$ and $H_2 (R,R) = \Omega^2_R$. The class $\tau_R$
is represented by the 2-cocycle
$$c\colon\ R\otimes R  \rightarrow \Omega^2_R,
\quad  c(f,g)=df\wedge dg.$$
The universal central extension $R^{\tau}$ is thus explicitly realized as
the algebra $R \oplus \Omega^2_R$ with multiplication
$$(f_1,\omega_1) (f_2,\omega_2)= (f_1f_2, f_1\omega_2 + f_2\omega_1 + df_1
\wedge df_2).$$

\vskip .2cm

\noindent
{\bf (1.4) Smooth and $d$-smooth algebras.} Let ${\cal A}lg$ be
the category of all associative algebras and ${\cal C}om$ be the subcategory of
commutative algebras. Let also ${\cal N}_d$ be the category of NC-nilpotent
algebras of
degree $d$ and ${\cal N}=\bigcup\limits_{d\ge 0} {\cal N}_d$. Then we have the
inclusions:
$${\cal C}om = {\cal N}_0 \subset {\cal N}_1 \subset
\cdots \subset {\cal N}\subset {\cal A}lg.\leqno (1.4.1)$$

\proclaim (1.4.2) Definition. Let $\cal C$ be one of the
categories in (1.4.1). A covariant functor $h\colon{\cal C}  \rightarrow
{\cal S}ets$ is called formally smooth, if the following equivalent
conditions hold: \hfill\break
{(i)} For any surjection $p\colon\ \Lambda^{\prime}  \rightarrow
\Lambda$ in $\cal C$ whose kernel is a nilpotent ideal, the map $h(p)\colon\
h(\Lambda^{\prime})  \rightarrow h(\Lambda)$ is surjective.\hfill\break
{(ii)} For any Abelian extension $p\colon\ \Lambda^{\prime}  \rightarrow
\Lambda$ in $\cal C$ the map $h(p)$ is surjective.

\proclaim (1.4.3) Proposition. Let ${\cal C} = {\cal N}$ or ${\cal
N}_d$
for some $d$. A functor $h\colon\ {\cal C}  \rightarrow {\cal S}ets$ is
formally smooth if an only if $h(p)$ is surjective for any central
extension $p\colon\ 
\Lambda^{\prime}  \rightarrow \Lambda$ in $\cal C$.

\noindent{\sl Proof:} See Proposition 1.2.3(b).

\vskip .1cm

For every associative algebra $R$ we have the functor
$$h^R\colon {\cal A}lg  \rightarrow {\cal S}ets,\qquad \Lambda \mapsto
\hbox{Hom}_{ {\cal A}lg} (R,\Lambda),$$
represented by $\Lambda$. If $\cal C$ is one of the categories above, let $h^R_{\cal C}$
be the restriction of $h^R$ to $\cal C$.
The following concept was studied in [Sche] [CQ].

\proclaim (1.4.4) Definition. An associative algebra $R$ is
called smooth (or quasi-free), if $h^R$ is formally smooth.

We will consider the following version of this concept.

\proclaim (1.4.5) Definition. Let $R$ be an associative
algebra. \hfill\break
{(a)} For $d\ge 0$ we say that $R$ is $d$-smooth, if $R$ is finitely
generated, $R\in {\cal N}_d$ and $h^R_{{\cal N}_d}$ is formally smooth.
\hfill\break
{(b)} We say that $R$ is NC-smooth, if $R$ is NC-complete and
$R/F^{d+1}R$ is $d$-smooth for each $d$.

For example, a 0-smooth algebra is the same as a smooth commutative
algebra, i.e. the coordinate algebra of a smooth affine algebraic variety
(Grothendieck's criterion). If $R$ is a $d$-smooth or NC-smooth algebra, we
denote by $\hbox{dim}(R)$ the dimension of the algebraic variety
$\hbox{Spec}(R_{ab})$. It is clear that for a $d$-smooth $R$ each
$R/F^{i+1}R, i\le d$ is $i$-smooth. Given a smooth commutative algebra $A$,
we call a $d$-smooth thickening of $A$ a $d$-smooth algebra $R$ together
with an isomorphism $R_{ab}  \rightarrow A$.

\proclaim (1.4.6) Proposition. Let $R$ be a finitely
generated quasi-free algebra. Then $R/F^{d+1}R$ is $d$-smooth for any
$d\ge 0$, and $R_{\ldb ab\rdb}$ is NC-smooth.

\noindent {\sl Proof:}  Follows  from $h^R_{{\cal
N}_d}=h^{R/F^{d+1}}_{{\cal N}_d}$.

\vskip .1cm

Given two algebras $R,S$, we denote by $R*S$ their free product, and by
$R\mathop{{}\hat{*}{}}S =(R*S)_{\ldb ab\rdb}$ its NC-completion.

\proclaim (1.4.7) Proposition. (a) If $R,S$ are $d$-smooth,
then so is $(R*S)/F^{d+1}(R*S)$. \hfill\break
(b) If $R,S$ are NC-smooth, then so is $R \mathop{{}\hat{*}{}}S$.

{\bf\it Proof.}\enspace\enspace  (a) $(R*S)/F^{d+1}(R*S)$ is the
categorical product of $R$ and $S$ in ${\cal N}_d$. Thus the functor
$h^{(R*S)/F^{d+1}}_{{\cal N}_d}$ takes $\Lambda \mapsto
h^R_{N-d}(\Lambda)\times h^S_{{\cal N}_d}(\Lambda)$ and it is formally smooth.

(b) Follows form (a) and from the identification
$$(R*S)/F^{d+1}(R*S)\simeq
\left((R/R^{d+1}R)*(S/F^{d+1}S)\right)\big/F^{d+1}.$$

\vskip .1cm

\proclaim (1.4.8) Corollary. If $R_1, \ldots, R_n$ are
coordinate algebras of smooth affine algebraic curves, then
$R_1\mathop{{}\hat{*}{}}\cdots\mathop{{}\hat{*}{}}R_n$ is NC-smooth.

\noindent {\sl Proof:} Follows from the fact that each $R_i$ is
quasi-free.

\vskip .3cm

\noindent
{\bf (1.5) Completions of $d$-smooth algebras.} Let ${\bb
C} {\<}x_1,\ldots, x_n{\>}  = {\bb C}[x_1]*\cdots* {\bb C}[x_n]$ be the
free
associative algebra on generators $x_1,\ldots, x_n$. We can view its
elements as noncommutative polynomials. Denote by ${\bf m}\subset {\bb C}
{\<}x_1,\ldots,x_n{\>}$ the two-sided ideal generated by the $x_i$. The
$\bf m$-adic
completion of ${\bb C}{\<}x_1, \ldots, x_n{\>}$ will be denoted ${\bb
C}{\<\!\<}x_1,\ldots, x_n{\>\!\>}$ and called the ring of noncommutative power
series.
This is a local ring whose maximal ideal $\hat{\bf m}$ is generated by the
$x_i$.

Suppose we have an associative algebra $R$, and $x\in \hbox{Spec} (R_{ab})$
is a ${\bb C}$-point. Let ${\bf m}_{x,ab}\subset R_{ab}$ be the corresponding
maximal ideal, and ${\bf m}_x = \pi^{-1}(m_{x,ab})\subset R$ be its preimage. The
${\bf m}_x$-adic completion of $R$ will be denote by $\hat{R}_x$. Its maximal
ideal will be denoted by $\hat{\bf m}_x$.

\proclaim (1.5.1) Proposition. (a) Let $R$ be $d$-smooth,
$\hbox{dim}(R)=n$, and $x\in \hbox{Spec}(R_{ab})({\bb C})$. Then $\hat{R}_x
\simeq
{\bb C}{\<\!\<} x_1,\dots, x_n{\>\!\>}/F^{d+1}$.\hfill\break
(b) Let $R$ be NC-smooth, $\hbox{dim}R=n$ and $x$ be as above.
Then $\hat{R}_x \simeq {\bb C}{\<\!\<}x_1,\ldots, x_n{\>\!\>}$.

\noindent 
{\sl Proof:} (a) ${\bb C}{\<}x_1,\ldots, x_n{\>}/F^{d+1}$
is
the free algebra in ${\cal N}_d$ generated by $x_x,\dots,x_n$, and ${\bb C}
{\<\!\<}
x_1,\ldots, x_n{\>\!\>}/F^{d+1}$ is its completion. If $\hbox{dim}(R) =n$, then
${\bf m}_x/{\bf m}^2_x$ is an $n$-dimensional vector space. Our statement is proved in
the same way as the formal tubular neighborhood theorem of Cuntz and
Quillen ([CQ], \S 6, Th. 2). In fact, it is a version of that theorem but
for algebras satisfying the polynomial identities expressed by $F^{d+1}=0$.
Part (b) follows from (a).

\proclaim (1.5.2) Corollary. Let $R,x,n$ be as in (1.5.1).
Then:\hfill\break
{(a)} If R is $d$-smooth and $i\le d$, then $R/{\bf m}_x^{i+1} \simeq {\bb
C} {\<}x_1,\ldots, x_n{\>}/{\bf m}^{i+1}$.\hfill\break
\item{(b)} If $R$ is NC-smooth and $i\ge 0$, then
$$R/{\bf m}_x^{i+1} \simeq {\bb C} {\<}x_1,\ldots, x_n{\>}/{\bf m}^{i+1}.$$
In any of these cases we have a natural isomorphism ${\bf m}_x^i/{\bf m}_x^{i+1} \simeq
({\bf m}_x/{\bf m}_x^2)^{\otimes i}$.

\vskip .2cm

\noindent
{\bf (1.6) Existence and uniqueness of thickenings.} We now formulate
the main result of this section.

\proclaim (1.6.1) Theorem. Let $A$ be a finitely generated
smooth commutative algebra. Then:\hfill\break
{(a)} For any $d\ge 0$ $A$ possesses a $d$-smooth thickening $R
\rightarrow A$ (with $R_ {ab}  \rightarrow A$ being an isomorphism). This
thickening is unique up to an isomorphism identical on $A$.\hfill\break
{(b)} Similar statement for NC-smooth thickenings.

\noindent {\sl Proof:} Only (a) needs to be proved, since (b) is a
formal consequence. So we start with the proof of uniqueness in (a). Let
$R^{\prime}, R^{\prime\prime}$ be  two $d$-smooth thickenings of $A$. By
induction we may assume that $R^{\prime}/F^dR^{\prime}$ and
$R^{\prime\prime}/F^dR^{\prime\prime}$ are isomorphic and identify them
both with the same $(d-1)$-smooth thickening $R_{d-1}$, as in the following
diagram:

$$\vbox{\setbox1=\hbox{$I^{\prime\prime}$}
\setbox2=\hbox{$f{}{\downarrow\uparrow}{}g$}
\setbox3=\hbox{$I_d{}{\downarrow}{}I_d$}
\halign{\tabskip 0pt
${}#{}$&${}#{}$&\hbox to \wd1{\hfil${}#{}$\hfil}&
${}#{}$&$\hbox to \wd2{\hfil${}#{}$\hfil}$&${}#{}$&
$\hbox to \wd3{\hfil${}#{}$\hfil}$&${}#{}$&${}#{}$\tabskip0pt\cr
0&\rightarrow&I^{\prime}&\rightarrow&R^{\prime}&\rightarrow&R_{d-1}&\rightarrow&0\cr
&&{\downarrow\uparrow}&&f{}{\downarrow\uparrow}{}g&&\phantom{\rm      Id}{}{\downarrow}{}{\rm Id}&&\cr
0&\rightarrow&I^{\prime\prime}&\rightarrow&R^{\prime\prime}&\rightarrow&R_{d
-1}&\rightarrow&0\cr}}$$

\noindent
Because both $R^{\prime},R^{\prime\prime}$ are $d$-smooth there exist
morphisms of extensions $f,g$ as shown. Since $R^{\prime}_{ab} =
R^{\prime\prime}_{ab} = (R_{d-1})_{ab}$, we find that $f$ and $g$ are isomorphisms by
Corollary 1.2.7.

We now prove the existence. It is  enough to show the following more
precise fact.

\proclaim (1.6.2) Proposition. Let $R$ be a $d$-smooth
algebra and $R^{\tau}$ be its universal central extension. Then
$R^{\tau}$ is $(d+1)$-smooth.

It is enough to prove that any surjection $\pi_{d+1}\colon S_{d+1}
\rightarrow R^{\tau}$ with nilpotent kernel and $S_{d+1} \in N_{d+1}$,
splits, i.e., admits $\sigma_{d+1}\colon\ R^{\tau}  \rightarrow S_{d+1}$
with $\pi_{d+1}\sigma_{d+1}=\hbox{Id}$. Let $S_d = S_{d+1}/F^{d+1}S_{d+1}$
and $q\colon\ S_{d+1}  \rightarrow S_d$ be the natural projection. Then
$S_d\in {\cal N}_d$. We will now gradually construct the following diagram:

$$\vbox{
\halign{\tabskip 0pt
$#$&
$#$&
$#$&
$#$&
$#$&
$#$&
$#$\tabskip 0pt\cr
&&\multispan{3}{$\displaystyle\hfil\mathop{{}\longrightarrow{}}^{\pi_{d+1}}
\hfil$}\cr
\noalign{\vskip -13pt}
&S_{d+1}&&&&R^{\tau}\cr
&&\displaystyle\mathop{{\nwarrow}}^{\alpha}&&\displaystyle\mathop{\swarrow}^
{\gamma}\cr
q&{}\big\downarrow{}&&U&&{}\big\downarrow{}&p\cr
&&&&\displaystyle\mathop{\searrow}^{\beta}\cr
&S_{d\phantom{{}+1}}&&&&R^{\phantom{\tau}}\cr
\noalign{\vskip -13pt}
&&\multispan{3}{$\displaystyle\hfil\mathop{{}\longleftarrow{}}^{\sigma_d}\hfil$}\cr
\noalign{\vskip -7pt}
&&\multispan{3}{$\displaystyle\hfil\mathop{{}\longrightarrow{}}_{\pi_d}\hfil
$}\cr}}$$

\noindent
Here $\pi_d$ is the surjection induced by $\pi_d$ (it exists becaue $R\in
{\cal N}_d$). Since $R$ is $d$-smoooth, we have $\sigma_d$ with
$\pi_d\sigma_d =
\hbox{Id}$. Let $U$ be the fiber product of $S_{d+1}$ and  $R$ over $S_d$
with
respect to $q$, $\sigma_d$, and $\alpha,\beta$ its natural projections. Let
$I=\hbox{Ker}(q)$ and $H=H_2(R,R_{ab})=\hbox{Ker}(p)$.
Then $\beta\colon\ U\rightarrow R$ is a central extension with kernel $I$
(since $q\colon\ S_{d+1}\rightarrow S_d$ is). Further, $\pi_{d+1}\alpha
\colon\ U\rightarrow R^{\tau}$ is a morphism of central extensions because
$p \pi_{d+1} \alpha=\pi_dq\alpha=\pi_d\sigma_d\beta=\beta$. Now, by
Proposition 1.3.8, there exist a morphism of extensions $j\colon\
R^{\tau}  \rightarrow U$. Define  $\rho_{d+1} = \alpha j\colon\ R^{\tau}
\rightarrow S_{d+1}$.

\proclaim (1.6.3) Lemma. $\pi_{d+1}\rho_{d+1}\colon\
R^{\tau}  \rightarrow R^{\tau}$ is an isomorphism.

Given the lemma, in order to construct $\sigma_{d+1}$, we need just to
compose $\rho_{d+1}$ with the inverse of $\pi_{d+1}\rho_{d+1}$.

\vskip .1cm

\noindent {\sl Proof of the lemma:}  First, notice that the
endomorphism of $R$ induced by $\pi_{d+1}\rho_{d+1}$, is the identity, i.e.
$p\pi_{d+1}\rho_{d+1} =p$. Indeed,

$$p\pi_{d+1} \rho_{d+1} = p\pi_{d+1}\alpha\gamma=\pi_d q \alpha \gamma =
\pi_d\sigma_d\beta \gamma = \beta \gamma =p.$$

\noindent
Further, the endomorphism of $H=\hbox{Ker}(p)$ induced by
$\pi_{d+1}\rho_{d+1}$ is the identity by Proposition 1.3.8. So
$\pi_{d+1}\rho_{d+1}$ is an isomorphism.

This completes the proof of Proposition 1.6.2 and Theorem 1.6.1.

\noindent
{\bf (1.6.4) Remarks.} Thus for a smooth finitely
generated commutative algebra $A$ the unique $d$-smooth thickening can be
constructed by repeatedly taking the universal central extension, starting
from $A$. One  can compare this with a similar situation in group theory.

If $G$ is a group whose Abelianization $G_{ab}=H_1(G,{\bb Z})$ is free, then
the universal coefficients formula gives a tautological class $\tau_G\in
H^2(G, H_2(G,{\bb Z}))$ and we can form the ``universal''central extension
with kernel $H_2(G,{\bb Z})$. Taking a free Abelian group $A$ with $n$
generators and  applying this construction $d$ times, we get the free
$(d+1)$-stage nilpotent group on  $n$ generators. See [Ev][BE] where this
folklore fact is implicit. As in (1.3.9), the very first universal central
extension (with kernel $\bigwedge^2A)$ is functorial in $A$, but the subsequent
ones are not.

\vfill\eject

\centerline {\bf \S2. NC-schemes.}

\vskip 1cm

\noindent
{\bf (2.1) Localization of NC-nilpotent algebras.} We
start with the following obvious remark.

\proclaim (2.1.1) Proposition. If $R$ is an NC-complete
algebra, $A=\|a_{ij}\|\in \hbox{Mat}_m(R)$ be a square matrix such that
$A_{ab} = \|(a_{ij})_{ab}\|\in \hbox{Mat}_m(R_{ab})$ is invertible, then
$A$ is invertible.

\noindent {\sl Proof:} Let $B\in \hbox{Mat}_m(R)$ be such that
$B_{ab}=A^{-1}_{ab}$, i.e., $U=AB-1\in \hbox{Mat}_m(F^1R)$. Now we find
$A^{-1}=\sum_{d=0}^{\infty} B\cdot (-1)^d\cdot U^d$, where the series
converges in the topology on $\hbox{Mat}_m(R)$ induced by the NC-topology
on $R$.

\vskip .1cm

\proclaim (2.1.2) Corollary. If $R$ is NC-complete and
$R_{ab}$ is local, then $R$ is local.

\noindent We now recall the framework of Ore localization [Ste].

\proclaim (2.1.3) Definition.  Let $R$ be an associative
algebra. A multiplicative subset $S\subset R-\{0\}$ is said to satisfy the
Ore conditions, if the following hold:\hfill\break
({OL1}) For any $a\in R, s\in S$ there are $b\in R, u\in S$ such that
$ua=bs$.\hfill\break
({OL2}) If $as=0$ with $a\in R,s\in S$, then $ta=0$ for some $t\in S$.
\hfill\break
({OR1}) For any $b\in R, u\in S$ there are $a\in R,s\in  S$ such that
$ua=bs$.\hfill\break
({OR2}) If $sa=0$ with $a\in R, s\in S$, then $at=0$ for some $t\in S$.

If $S$ satisfies the Ore conditions one defines the algebra of fractions
$R[S^{-1}]$ as the set of equivalence classes of, say, right fractions
$u^{-1}b$, see [Ste], Ch. II. This algebra has the following properties.

\proclaim (2.1.4) Proposition. (a) There exist a natural
homomorphism $c\colon\ R  \rightarrow R[S^{-1}]$ taking any $s\in S$ into an
invertible element.\hfill\break
(b) For any algebra homomorphism $\varphi\colon\ R  \rightarrow R^{\prime}$
such that $\varphi(s)$ is invertible for any $s\in S$, there is a unique
homomorphism $\overline{\varphi}\colon\ R[S^{-1}]  \rightarrow R$ such that
$\varphi =\overline{\varphi}c$. \hfill\break
(c) $R[S^{-1}]$ is flat as a left and as a right $R$-module.
\hfill\break
(d) $c(a)=0$ if and only if $as=0$ for some $s\in S$.

All these properties can be found in [Ste], Ch. II. More precisely, (c) is
Prop. 3.5 of {\it loc. cit.}, (b) is Prop. 1.1 and (a),(d) are properties
F1,F3 which define the abstract concept of a ring  of fractions (of which
the Ore construction proves the existence).

\vskip .1cm

\proclaim (2.1.5) Proposition. Let $R$ be NC-nilpotent,
$\pi\colon\ R  \rightarrow R_{ab}$ the Abelianization map and
$\overline{S}\subset R_{ab} -\{0\}$ be any multiplicative subset. Then
$S=\pi^{-1}(\overline{S})$ satisfies the Ore conditions.

\noindent {\sl Proof:} (OL1) For every $n\ge 0$ we have the
identity
 $$
s^{n+1}a=\left(\sum^{n}_{i=0}s^{n-i}ad(s)^i(a)\right)s+ad(s)^{n+1}(a).
\leqno (2.1.5.1)$$
If $F^{n+1}R=0$, the last term on the right vanishes. So taking
$$u=s^{n+1}, b = \sum^n_{i=0}s^{n-i}ad(s)^i(a),$$
we get $ua=bs$.

(OL2) If $as=0$, then we find inductively:
$$sa=[s,a], s^2a=[s,sa]=[s,[s,a]],\ldots,s^na=ad(s)^na.$$
Thus if $F^nR=0$, then $s^na=0$.
The proof of (OR1-2) is similar.

\vskip .1cm

\proclaim (2.1.6) Theorem. Let $R$ be NC-nilpotent,
$\overline{S}\subset R_{ab}-\{0\}$ a multiplicative subset, and
$S=\pi^{-1}(\overline{S})\subset R$ its preimage. Then $R[S^{-1}]$ is
NC-nilpotent and
$$\hbox{gr}^d_F(R[S^{-1}])=\hbox{gr}^d_F(R)[\overline{S}^{-1}],$$
where on the right we have the usual module of fractions of the
$R_{ab}$-module $\hbox{gr}^d_F(R)$.

The proof is based on the following fact.

\vskip .1cm

\proclaim (2.1.7) Proposition. We have
$F^d(R[S^{-1}])=S^{-1}\cdot (F^dR)$.

We first deduce the theorem from the proposition. Because of the properties
of modules of fractions with respect to Ore sets ([Ste], Ch. II. \S3), we can
write $F^d(R[S^{-1}])= R[S^{-1}]\mathop{{}\otimes{}}\limits_{R}(F^dR)$.
Further,
since $R[S^{-1}]$ is flat over $R$ by (2.1.4)(c), we have:
$$\eqalign{
\hbox{gr}^d_F(R[S^{-1}])&= R[S^{-1}]\mathop{{}\otimes{}}_{R}\hbox{gr}^d_F
(R) =
(R[S^{-1}]\mathop{{}\otimes{}}_{R}R_{ab})\mathop{{}\otimes{}}_{R_{ab}}\hbox{
gr}^d_F(R)=\cr
& =
R_{ab}[\overline{S}^{-1}]\mathop{{}\otimes{}}_{R_{ab}}\hbox{gr}^d_F(R)=\hbox
{gr}^d_F(R)[\overline{S}^{-1}]\cr}.$$

We now prove Proposition 2.1.7. It is obviously enough to prove that
$$ [s^{-1}_0a_o, [s^{-1}_1 a_1,\ldots, [s^{-1}_{d-1} a_{d-1}, s^{-1}_d
a_d]\ldots]\in S^{-1} F^dR$$
for any $s_i\in S,a_i\in R$. But this is proved by induction, using the
identities
$$[xy,z]=x[y,z]+[x,z]y,\qquad
[s^{-1},z]=-s^{-1}[s,z]s^{-1}, \leqno (2.1.7.1)$$
$$as^{-1}=s^{-1}a+s^{-2}[s,a]+s^{-3}[s,[s,a]]+\cdots\leqno (2.1.7.2)$$
of which the third one is obtained by iterating the second one (or from
(2.1.5.1)). Theorem is proved.

\vskip .1cm

In the sequel we will write $R[\overline{S}^{-1}]$ for the localization
$R[\pi^{-1}(\overline{S})^{-1}]$, where $R$ is a NC-nilpotent algebra and
$\overline{S}\subset R_{ab}$ a multiplicative subset. In particular, we
write $R[g^{-1}]$ for $R[\pi^{-1}\{g^i\})^{-1}]$, $g\in R_{ab}$.

\vskip .1cm

\proclaim (2.1.8) Definition. For an NC-complete algebra $R$
and a multiplicative subset $T\subset R_{ab}$ we define the localization
$$R[\![ T^{-1}]\!] = \lim_{\leftarrow} (R/F^{d+1}R)[T^{-1}].$$

Because localization  does not commute with inverse limits, not all
elements of $R[\![T^{-1}]\!]$ can be represented as actual fractions.

\vskip .3cm

\noindent
{\bf (2.2) NC-schemes.} Let $R$ be an NC-nilpotent
algebra.
Denote by $X_{ab} =\hbox{Spec}(R_{ab})$ the affine scheme
corresponding to the commutative algebra $R_{ab}$.
We will now  construct a ringed space $X=\hbox{Spec}(R)=(X_{ab}, {\cal
O}_X)$
with underlying space $X_{ab}$, i.e., the set of prime ideals in $R_{ab}$.
As well known, the basis of topology in $X_{ab}$ is formed by the principal
open subsets $D_g = \{\wp \in \hbox{Spec}(R_{ab})\colon\
g\notin \wp\}$ for $g\in R_{ab}$. We define a presheaf $\tilde{{\cal O}}$ on
this basis of topology by putting $\tilde{{\cal O}}(D_g)=R[g^{-1}]$, the
localization defined in (2.1). Then we set ${\cal O}_X$ to be the associated
sheaf.
Explicitly, this means that we first form the stalks
$${\cal O}_\wp\colon\ = R_\wp = \mathop{\lim_{\rightarrow}}_{\wp\in D_g}
R[g^{-1}],\qquad \wp\in Spec(R_{ab}),$$
then make $\coprod {\cal O}_\wp$ into a covering space of $X_{ab}$ in a
standard way
and define ${\cal O}_X$ as the sheaf of continuous sections.

\vskip .1cm

\proclaim
 (2.2.1) Proposition. (a) Each stalk ${\cal O}_p$ is
a
local ring. \hfill\break
(b) We have $\Gamma(X_{ab}, {\cal O}_x)=R$.

\noindent {\sl Proof:} (a) Follows from (2.1.2).

 (b) We are reduced to the following situation. Given a covering
$X_{ab}=\bigcup\limits_{g\in J}D_g$, it is required to prove that the \v
Cech complex
$$R  \rightarrow \prod_{g\in J}R[g^{-1}]\rightarrow \prod_{g_1,g_2\in J}
R[g^{-1}_1,g^{-1}_2]$$
is exact in the middle term. The NC-filtration makes it into a complex of
filtered vector spaces with finite filtrations. The associated graded
complex is a direct sum of complexes of the form
$$M  \rightarrow \prod_{g\in J} M[g^{-1}]  \rightarrow \prod_{g_1,g_2\in
J} M[g^{-1}_1, g^{-1}_2]$$
where $M=\hbox{gr}^d_F(R)$ is an $R_{ab}$-module. The exactness of such a
complex is well known. So the original complex is exact as well, proving
out statement.

\vskip .1cm

\proclaim (2.2.2) Definition. Let $R$ be an NC-complete
algebra, $X_{ab}=\hbox{Spec}(R_{ab})$. The formal spectrum $\hbox{Spf}(R)$
is the ringed space $(X_{ab},{\cal O}_X)$ where ${\cal O}_X$ is the sheaf of
topological
rings obtained as the inverse limit of the structure sheaves of
$\hbox{Spec}(R/F^{d+1}R)$.

Ringed spaces of the form $\hbox{Spf}(R)$ will be called affine NC-schemes.
Let ${\cal NC}$ be the category of NC-complete algebras.

\proclaim (2.2.3) Proposition. The functor $X\mapsto \Gamma
(X_{ab},{\cal O}_X)$ establishes an equivalence between the category of
affine
NC-schemes and the opposite category to $\cal NC$.

\noindent {\sl Proof:} Follows from (2.2.1).

\vskip .1cm

\proclaim  (2.2.4) Definition. An NC-scheme is a ringed space
$X=(M,{\cal O}_X)$ locally isomorphic to an affine NC-scheme.

Thus every NC-scheme $X$ gives rise to an ordinary scheme
$X_{ab}=(M,{\cal O}_{X,ab})$. Nontrivial examples of NC-schemes will be
given in
\S5.

\vskip .1cm

\proclaim (2.2.5) Definition. An NC-scheme $X$ is called
NC-nilpotent (of degree $d$) if ${\cal O}_X$ is. We say that $X$ is of
finite
type, if it is NC-nilpotent, $X_{ab}$ is a scheme of finite type over $\bb
C$ and the sheaves $\hbox{gr}^d_F{\cal O}_X$ on $X_{ab}$ are coherent.

Recall that ${\cal C}om\subset {\cal N}$ denote the categories of
 commutative and
NC-nilpotent algebras. An NC-scheme $X$ defines a covariant functor
 $$\tilde{h}_X\colon\ {\cal NC}  \rightarrow
{\cal S}ets,\quad \tilde{h}_X(\Lambda)=\hbox{Hom}_{\rm NC-sch}
(\hbox{Spf}(\Lambda), X). \leqno (2.2.6)$$
Let $h_X$ be the restriction of $\tilde{h}_X$ to $\cal N$.
\vskip .1cm

\proclaim (2.2.7) Proposition. (a) $\tilde{h}_X$ commutes
with inverse limits in $\cal NC$. \hfill\break
(b) The restriction of $h_x$ to ${\cal C}om$
 is the functor represented by $X_{ab}$.\hfill\break
(c) The correspondence $X\mapsto h_X$ embeds the category of NC-schemes as a
full subcategory into the category of covariant functors ${\cal N}  \rightarrow
{\cal S}ets$.

\noindent {\sl Proof:} (a) The category ${\cal NC}^{op}$ being a full
subcategory of the category of NC-schemes, the statement follows from the
general fact that a representable contravariant functor commutes with
direct limits. Part (b) is clear. Part (c) follows from (a) and the fact
that a NC-complete algebra is an inverse limit of NC-nilpotent ones.

\vskip .3cm

\noindent
{\bf (2.3) Smooth and $d$-smooth NC-schemes.} Let
$X=(M, {\cal O}_X)$ be an NC-scheme. Denote $X^{\le d} = (M, {\cal
O}_X/F^{d+1})$ the
$d$th truncation  of $X$. Recall that ${\cal N}_d$ denotes the category of
NC-nilpotent algebras of degree $d$. Let $h_X^{{\cal N}_d}$ be the
restriction of
the functor  $h_X$ to ${\cal N}_d$.

\proclaim (2.3.1) Definition. (a) An NC-scheme $X$ is called
$d$-smooth if $X$ is of finite type, $F^{d+1}{\cal O}_X=0$ and the functor
$h_X^{\cal NC}$
is formally smooth (1.4.2). \hfill\break
(b) $X$ is called smooth, if $X^{\le d}$ is $d$-smooth for each $d$.

The definitions imply that the underlying scheme $M=X_{ab}$ of a
($d$-)smooth NC-scheme $X$ is a smooth algebraic variety. Given a smooth
algebraic variety M, we define the category ${\rm Th}^d(M)$ (resp.
${\rm Th}^{\infty}(M)$)
whose objects are $d$-smooth (resp. smooth) NC-schemes $X$ together with an
isomorphism $X_{ab}  \rightarrow { M}$, and morphisms are morphisms of
NC-schemes identical on ${ M}$. Objects of this category will be called
$d$-smooth (resp. smooth) thickenings of $M$.

\vskip .1cm

\proclaim (2.3.2) Proposition. Any morphism of
$\hbox{Th}^d(M)$, $d\le \infty$, is an isomorphism.

\noindent {\sl Proof:} Follows from Corollary 1.2.7 and Theorem
1.6.1.

\vskip .1cm

Let a covariant functor $h\colon\ {\cal N}_d  \rightarrow {\cal S}ets$
be given.  For any Cartesian diagram of algebras in ${\cal N}_d$
$$\vbox{
\halign{\tabskip 0pt
$#$&
$#$&
$#$&
$#$&
$#$\tabskip 0pt\cr
&\Lambda_{12}&\displaystyle\mathop{{}\longrightarrow{}}^{q_1}&\Lambda_{1}\cr
q_2&{\downarrow}&&{\downarrow}&p_1\cr
&\Lambda_{2\phantom{1}}&\displaystyle\mathop{{}\longrightarrow{}}_{p_2}&
\Lambda_{\phantom{1}}\cr
}}\leqno (2.3.3)$$
(i.e., $\Lambda_{12}\simeq \Lambda_1\mathop{{}\times{}}\limits_{\Lambda}
\Lambda_2$)
we have a natural map
$$j\colon\ h(\Lambda_{12})  \rightarrow
h(\Lambda_1) \mathop{{}\times{}}\limits_{h(\Lambda)} h(\Lambda_2).
\leqno (2.3.4)$$
The following theorem (cf. [Schl]) will be used to construct many examples
of thickenings.

\vskip .1cm

\proclaim (2.3.5) Theorem. Let $M$ be a smooth algebraic
variety and $h^{{\cal C}om}_M\colon\ {\cal C}om  \rightarrow {\cal S}ets$ the
corresponding representable functor. Then the category $\hbox{Th}^d(M)$
(resp. $\hbox{Th}^{\infty}(M)$) is equivalent to the category of formally
smooth functors $h\colon\ {\cal N}_d  \rightarrow {\cal S}ets$ (resp. ${\cal
N} \rightarrow {\cal S}ets$) such that $h\big\vert_{ {\cal C}om}= 
h^{{\cal C}om}_M$
and satisfying the following left exactness properties: \hfill\break
{(1)} If, in (2.3.3), $\Lambda_1=\Lambda_2$ and $p_1=p_2$ is a central
extension, then $j$ is a bijection.\hfill\break
{(2)} If $\Lambda$ is commutative and $\Lambda_2 = \Lambda \oplus I$,
where $I$ is a $\Lambda$-module, see (1.2.4), then $j$ is a bijection.

\noindent {\sl Proof:} It is enough to consider $d<\infty$
so we assume this. What we
really need to do is to prove that any functor satisfying the listed
properties, is representable by an NC-scheme, nilpotent of degree $d$.
Further, it suffices to prove this when $M$ is affine. Indeed, the affine
case being established and given arbitrary
$M$, $h\colon\ {\cal N}_d  \rightarrow
{\cal S}ets$ satisfying the conditions of the theorem, we get a thickening
${ U}^{(d)}$ for any affine open ${ U}\subset M$, and it
 is defined up to a unique
isomorphism. This allows us to glue the $U^{(d)}$ together.

So we assume that ${ M}=\hbox{Spec}(A)$ is affine. We know that $M$ has
a
unique, up to an isomorphism, $d$-smooth thickening $X=\hbox{Spec}(R)$, see
Theorem 1.6.1, but we need to identify $h$ with the representable functor
$h^R\colon\Lambda \mapsto \hbox{Hom}(R,\Lambda)$.

Since $h$ is formally smooth and $\pi\colon\ R  \rightarrow A=R_{ab}$ is a
surjection  with nilpotent kernel, the identity $\xi_0\in h(A) = 
{\rm Hom}(A,A)$,
lifts to some $\xi\in h(R)$ which
gives a natural transformation $\xi_*\colon\ h^R  \rightarrow h$. We are
reduced to the following.

\proclaim (2.3.6) Lemma.  $\xi_*$ is an isomorphism of
functors, i.e., the map
$$\xi^{\Lambda}_*\colon\ \hbox{Hom}(R,\Lambda)  \rightarrow h(\Lambda),\qquad
{ f} \mapsto h(f)(\xi),$$
is a bijection for any $\Lambda \in {\cal N}_d$.

\noindent {\sl Proof:} We  use the induction the degree of
NC-nilpotency of $\Lambda$. So we assume that the statement is true for any
$\Lambda
\in {\cal N}_{d-1}$ and consider a central extension
$$0 \rightarrow I  \rightarrow \Lambda^{\prime} \mathop{{}\rightarrow{}}^{p}
\Lambda  \rightarrow 0,\qquad \Lambda \in {\cal N}_{d-1}.$$
Let us prove the bijectivity of $\xi^{\Lambda^{\prime}}_*$. By  Proposition
(1.2.5)(b) and the condition (2) of our theorem, we have a map
 $$\displaystyle h(\Lambda^{\prime} )
\mathop{{}\times{}}_{h(\Lambda_{ab})} h(\Lambda_{ab} \oplus I)  \rightarrow
h(\Lambda^{\prime}) \mathop{{}\times{}}_{h(\Lambda)} h(\Lambda^{\prime}).
\leqno (2.3.7)$$
For every $\eta \in h(\Lambda)$ let $\eta_{ab}\in h(\Lambda_{ab})$ be the
image of $\eta$ under $\pi\colon\ \Lambda \rightarrow \Lambda_{ab}$. Note
that $\eta_{ab}$ is just a morphism $A  \rightarrow \Lambda_{ab}$.
Given $g\colon\ A  \rightarrow \Lambda_{ab}$, the preimage of $g$ under
$\hbox{Hom}(A, \Lambda_{ab}\oplus I)  \rightarrow \hbox{Hom}(A,
\Lambda_{ab})$, is $\hbox{Der}(A,I)$ (here the structure of A-module on I is
given by $g$). Thus the map (2.3.7) defines, for any $\eta \in h(\Lambda)$,
the action of $\hbox{Der}(A,I)$ (with A-module structure on $I$ given by
$\eta_{ab}$) on $h(p)^{-1}\eta \subset h(\Lambda^{\prime})$. Properties (1)
and (2) imply that this action makes $h(p)^{-1}\eta$ into a principal
homogeneous space over $\hbox{Der}(A,I)$. The same conclusion holds for
$h^R$. But since $\xi^{\Lambda^{\prime}}_*$ takes, for any $\eta \in
h(\Lambda)$,
$$h(p)^{-1} (\eta)  \rightarrow
h^{\Lambda^{\prime}}(p)^{-1}\left(\xi^{\Lambda}_*(p)\right)$$
and is a morphism of principal homogeneous spaces, it is a bijection.
Theorem is proved.

\vfill\eject

\centerline {\bf \S3. The affine NC-space and Feynman-Maslov operator
calculus.}

\vskip 1cm

\noindent
{\bf (3.1) Free associative algebras.}Let $V$ be a
vector space and ${\cal A}ss(V)= \mathop{\bigoplus}\limits_{d\ge 0}
V^{\otimes d}$ be
the tensor (free associative) algebra generated by $V$. Clearly,
${\cal A}ss (V)_{ab} = S(V)$ is the symmetric algebra of $V$. When
$V={\bb C}^n$
with basis $x_1,\ldots, x_n$, we identify ${\cal A}ss(V)$ with ${\bb
C}\<x_1,
\ldots, x_n\>$ and $S(V)$ with ${\bb C}[x_1,\ldots, x_n]$.

\vskip .1cm

\proclaim (3.1.1) Definition. The affine NC-space $A^n_{NC}$
is the NC-scheme $\hbox{Spf}({\bb C}{\<} x_1, \ldots, x_n{\>}_{\ldb
ab\rdb})$.

\vskip .1cm

Thus, as a ringed space, $A^n_{NC}$ is the ordinary affine space $A^n$
equipped  with a certain sheaf ${\cal O}^{NC}$ of noncommutative algebras
on the Zariski topology of $A^n$. The aim of this section is to describe
${\cal O}^{NC}$ explicitly by using the so-called calculus of ordered
operators as developed by Feynman, Maslov and Karasev [Fe][Mas][KM].
The essence of this approach is to represent elements of ${\bb
C}{\<}x_1,\ldots,x_n{\>}$ by ordinary polynomials in $N\ge n$ variables. Let
$[n]=\{1,2,\ldots,n\}$. Let $\tau\colon\ [N]  \rightarrow [n]$ be any map
and ${f} \in {\bb C}[{ y}_1,\ldots,y_N]$ be a polynomial. Define a
noncommutative polynomial
$$\ldb f(\mathop{x}^{\tau})\rdb = \ldb f({\mathop{x}^{1}}_{\tau(1)},\ldots,
{\mathop{x}^{N}}_{\tau(N)})\rdb \in {\bb C} {\<}x_1,\ldots,x_N{\>}$$
by replacing each monomial $c\cdot y_1^{i_1}\ldots y_N^{i_N}$ in ${f}$
by
$c\cdot x^{i_1}_{\tau(1)} x_{\tau(2)}^{i_2}\ldots x_{\tau(n)}^{i_n}$. For
example, let $N=n=2$ and ${f}(y_1,y_2)=(y_1+y_2)^2$. Taking $\tau=
\hbox{Id}$ and $\sigma = (12)$ we get
$$\ldb f(\mathop{x}^{\tau})\rdb= ({\mathop{x}^{1}}_1 +
{\mathop{x}^{2}}_2)^2=
x^2_1 + 2x_1x_2 + x^2_2,\quad
\ldb f(\mathop{x}^{\sigma})\rdb= ({\mathop{x}^{2}}_1 +
{\mathop{x}^{1}}_2)^2=
x^2_1 + 2x_2x_1 + x^2_2.$$
Thus the numbers over the variables indicate the order. Given $f\in {\bb C}
[x_1,\ldots, x_n]$ we will use the default notation for the standard
ordering ($\tau = \hbox{Id}$):
$$\ldb f(x)\rdb=\ldb f(x_1,\ldots, x_n)\rdb= \ldb f({\mathop{x}^1}_1,\ldots,
{\mathop{x}^{n}}_n)\rdb.$$

\proclaim
(3.1.2) Proposition. ${\bb C} {\<}x_1,\ldots,
x_n{\>}$
can be identified with the space of formal finite sums
$$\sum_{N,\tau} \ldb f_{\tau} (\mathop{x}^{\tau})\rdb,\quad f_{\tau}\in {\bb
C}[y_1,\ldots, y_N]$$
modulo the following cancellation rules: \hfill\break
{(C1)} If $f_{\tau}\in {\bb C}[y_1,\ldots, y_N]$ and
$\tau(i)=\tau(i+1)$, then
$$\ldb f_{\tau} (\mathop{x}^{\tau})\rdb = \ldb(r_if)(
{\mathop{x}^{1}}_{\tau(1)},\ldots,
{\mathop{x}^{i}}_{\tau(i)},
{\mathop{x}^{i+1}}_{\tau(i+2)},\ldots,
{\mathop{x}^{N-1}}_{\tau(N)})\rdb,$$
where $r_if\in {\bb C}[y_1,\ldots, y_{N-1}]$
is defined by
$$(r_if)(y_1,\ldots,y_{N-1})={f}(y_1,\ldots, y_i, y_i, y_{i+1},\ldots,
y_{N-1}).$$
{(C2)} If $g\in {\bb C}[y_1,\ldots, y_{N-1}]$ and $s_ig\in {\bb
C}[y_1, \ldots, y_N]$ is defined by
$$(s_ig)(y_1,\ldots, y_N)=g(y_1,\ldots,\hat y_i,\ldots, y_N),$$
then
$$\ldb (s_i g)(\mathop{x}^{\tau})\rdb = \ldb
g({\mathop{x}^{1}}_{\tau(1)},\ldots, {\mathop{x}^{i-1}}_{\tau(i+1)},
{\mathop{x}^{i}}_{\tau(i+1)},\ldots,
{\mathop{x}^{N-1}}_{\tau(N)})\rdb.$$

\noindent {\sl Proof:} Given any $n$ associative algebras,
$A_1,\ldots, A_n$, their free product $A_1 *\cdots *A_n$ is spanned by
formal products $a_{\tau(1)}*\cdots * a_{\tau(N)}$, where $N\in {\bb Z}_+$,
$\tau\colon\ [N]  \rightarrow[n]$, and $a_{\nu}\in A_{\nu}$.
These products are ${\bb C}$-multilinear and are subject to the two
cancellation rules:

\noindent {(C1$^{\prime}$)} $a_i*b_i = a_ib_i, a_ib_i\in A_i$.

\noindent {(C2$^{\prime}$)} $a_i*1_j*a_k = a_i*a_k$,  where
 $a_i\in A_i, a_k\in A_k$ and 
$1_j\in A_j$ is the unit.

So, $A_1*\cdots * A_n$ is a quotient of
$\mathop{\bigoplus}_{N,\tau}A_{\tau(1)}\otimes\cdots\otimes A_{\tau(N)}$. When
$A_1= \cdots = A_n = {\bb C}[x]$, we
have
$$A_1*\cdots *A_n={\bb C}{\<}x_1,\ldots,x_n{\>},\quad A_{\tau(1)}\otimes
\cdots
\otimes A_{\tau(N)}={\bb C}[y_1,\ldots, y_N]$$
and (C1$^{\prime}$), (C2$^{\prime}$) become (C1) and (C2).

\vskip .3cm

\noindent
{\bf (3.2) Free Lie algebras and the normal form in ${\bb
C}{\<}x_1,\ldots,x_n{\>}$.}  For a vector space $V$ let
${\cal L}ie (V)$ be the free
Lie algebra generated by $V$. It is graded:
$${\cal L}ie(V)=\mathop{\bigoplus}_{d\ge 1} {\cal L}ie ^d(V),\quad
{\cal L}ie^1(V)={ V}.$$
Every associative algebra, in particular, ${\cal A}ss (V)$, can be
regarded as a Lie
algebra via $[a,b]=ab-ba$. The following fact is classical.

\vskip .1cm

\proclaim (3.2.1) Theorem. (a) ${\cal L}ie(V)$ is
isomorphic
to the Lie subalgebra in ${\cal A}ss (V)$ generated by $V$.
\hfill\break
(b) The embedding ${\cal L}ie(V)\subset {\cal A}ss(V)$
identifies
${\cal A}ss (V)$ with the universal enveloping algebra  of
${\cal L}ie(V)$.

\proclaim (3.2.2) Corollary. Let $\{\beta_i\},
i=1,2,\ldots,$ be any ${\bb C}$-basis of ${\cal L}ie(V)$. Then the
ordered
monomials ${\beta_1}^{m_1}{\beta_2}^{m_2}\dots,$ where $m_i\in  {\bb Z}_+$
are 0 for almost all $i$, form a ${\bb C}$-basis in ${\cal A}ss(V)$.

\vskip .1cm

We now take $V={\bb C}^n$ with basis $x_1,\ldots, x_n$ and denote
${\cal L}ie(V)$ by ${\cal L}(x_1,\ldots, x_n)$ and ${\cal L}ie^d(V)$ by
${\cal L}^d(x_1,\ldots,x_n)$. Every element of ${\cal L}^d(x_1,\ldots,
x_n)$ can be represented (not uniquely) as a linear combination of Lie
monomials in $x_1,\ldots, x_n$.
For every $d\ge 2$ choose an ordered basis $B_d$ in ${\cal L}^d(x_1,\ldots,
x_n)$ (for example, one can take the Hall basis described in [MKS], \S5.6,
Ex. 10). Let $B=\coprod_{d\ge 2} B_d$ ordered lexicographically ($B_d$
precedes $B_{d^{\prime}}$ if $d < d^{\prime}$). For $\beta \in B_d$
we write $\hbox{ord}_{\rm NC} (\beta)=d-1$ (this is the total number of
bracket pairs). Let us number the elements of $B$ according to their total
order as $\beta_1,\beta_2,\ldots$. Denote by $\Lambda$ the set of all
functions $\lambda\colon\ B  \rightarrow {\bb Z}_+$ with finite support.
For $\lambda \in \Lambda$ we set
$$M_{\lambda}(x) =\beta_1^{\lambda(\beta_1)}\beta_2^{\lambda(\beta_2)}
\ldots \in {\bb C}{\<}x_1,\ldots, x_n{\>},\quad \hbox{ord}(\lambda)
=\sum_{\beta\in B} \lambda (\beta)\hbox{ord}_{\rm NC}(\beta)\in {\bb Z}.$$
The following is a reformulation of Corollary 3.2.2.

\vskip .1cm

\proclaim (3.2.3) Proposition. Every element of ${\bb C}
{\<}x_1,\ldots, x_n{\>}$ can be uniquely written as a finite sum
$$\sum_{\lambda \in \Lambda} \ldb f_{\lambda}(x)\rdb M_{\lambda}(x).$$

\vskip .2cm

\noindent
{\bf (3.3) Manipulations with ordered symbols.} Here we
recall several formulas due to Maslov and Karasev [Mas][KM].

For a polynomial $f\in {\bb C}[\xi]$ its difference derivative $\partial f$
is
defined by

$$\partial f(\xi^{\prime}, \xi^{\prime\prime})= {\partial f\over \partial
\xi}(\xi^{\prime},\xi^{\prime\prime})=
{f(\xi^{\prime})-f(\xi^{\prime\prime})\over
\xi^{\prime}-\xi^{\prime\prime}}.$$

The $m$-fold iterated difference derivative depends only on $m+1$
variables, reducing to
$${\partial^mf\over \partial \xi^m}(\xi^{(0)},\ldots,
\xi^{(m)})=\sum^m_{j=0}{f(\xi^{(j)})\over
\prod\limits_{i\not=j}(\xi^{(j)}-\xi^{(i)})}$$
Accordingly, for a polynomial  $f(\xi_1,\ldots, \xi_n)$ in $n$ variables
its partial difference derivative
$${\partial^{m_1+\cdots+m_n}f\over \partial \xi_1^{m_1}\cdots
\partial\xi_n^{m_n}}$$
depends on $(m_1+1)+\cdots + (m_n+1)$ variables $\xi_i^{(j)}, 0\le j\le
m_i$. Now the main formula of Maslov ([Mas], Ch. 0, Th. 4.3) is:
\vskip .1cm

\proclaim (3.3.1) Change of order formula I.  For any $g\in
{\bb C}[y_1,y_2]$ we have an equality in ${\bb C}{\<}a,b{\>}$:
$$\[ g(a^{\hskip -.15cm ^2}, b^{\hskip -.15cm ^1})\] - 
\[g(a^{\hskip -.15cm ^1}, b^{\hskip -.15cm ^2})\] = 
 \biggl[ \!\!\biggl[
[a,b]^{\hskip -.45cm ^3}
\hskip .4cm
{\delta^2 g\over \delta y_1 \delta y_2}(a^{\hskip -.15cm ^2}, 
a^{\hskip -.15cm ^4}, b^{\hskip -.15cm ^1}, b^{\hskip -.15cm ^5})
\biggr]\!\!\biggr].$$


\noindent
By iterating this, one gets a more precise formula involving the usual
partial derivatives ([KM] App. 1, Th. 1.9):

\proclaim (3.3.2) Change of order formula II. In the above
situation we have
$$\ldb g(\mathop{a}^{2}, \mathop{b}^{1})\rdb = \sum_{m,l\ge 0} \LDB
{\partial^{l+m}g\over \partial y_1^l\partial y_2^m}(\mathop{a}^{1},
\mathop{a}^{3}){\mathop{K}^{2}}_{l,m}\RDB,$$
where
$$K_{l,m} = {1\over l! m!} \[(a^{\hskip -.15cm ^1}- a^{\hskip -.15cm ^3})^l
(a^{\hskip -.15cm ^2} - b^{\hskip -.15cm ^4})^m\]$$
satisfies $\hbox{ord}_{\rm NC}
(K_{l,m}) \ge \max (l,m)$.

To write down $K_{l,m}$ one should first think of the
$a^{\hskip -.15cm ^i}$ as different variables, expanding
the powers by the binomial formula, and then make out
of every commutative monomial in $a^{\hskip -.15cm ^1}, 
a^{\hskip -.15cm ^2}, a^{\hskip -.15cm ^3}, b^{\hskip -.15cm ^4}$
a noncommutative monomial in $a,b$, as indicated by the superscripts.

 

One has similar formulas for interchanging the order of any two consecutive
variables in a many-variable symbol $g(y_1,\ldots, y_N)$, see e.g. [Mas].
For the next formula, see [KM], App. 1, n. (1.3).

\vskip .1cm

\proclaim (3.3.3) Commutation formula. For any $f\in {\bb C}
[y_1,\ldots, y_n]$ we have the  identity in ${\bb C}\< x_1,\ldots, x_n, a\>$:
$$\displaylines{a \cdot \ldb f(x_1, \ldots, x_n)\rdb =
\sum^{\infty}_{i_1,\ldots, i_n=0}
{1\over i_1!\ldots i_n!} \LDB {\partial^{i_1+\cdots + i_n}f\over \partial
y_1^{i_1}\ldots \partial y_n^{i_n}}(x_1,\ldots, x_n)\RDB\cdot\cr
{}\cdot {\rm ad}(x_n)^{i_n}\ldots {\rm ad}(x_1)^{i_1}(a).\cr}$$

\proclaim (3.3.4) Taylor formula. For any $f(y)\in {\bb
C}[y]$ we have the equality in ${\bb C}{\<}a,b{\>}$:
$$f(a+b) = \sum_{k\geq 0} \[f^{(k)}(a^{\hskip -.15cm ^1} + b^{\hskip -.15cm ^3})
X_k^{\hskip -.15cm ^2}\],$$
where 
$$X_k = {1\over k!} \[(\overline{a+b}^{\hskip -.5cm ^2}
 \hskip .45cm - a^{\hskip -.15cm ^1}
-b^{\hskip -.15cm ^3})^k \]\in {\bf C}\< a,b\>$$
satisfies $\hbox{ord}_{\rm NC} (X_k)\ge [(k+1)/2]$.


\vskip .2cm

\noindent
{\bf (3.4) Multiplication of elements in the normal
form.}
The formulas (3.3.1-3) can be used as an algorithm for bringing the product
$$\displaystyle \left(\sum_{\lambda\in
\Lambda} \ldb
f_{\lambda}(x)\rdb { M}_{\lambda}(x)\right)\left(\sum_{\mu\in
\Lambda}\ldb
g_{\mu} (x)\rdb M_{\mu} (x)\right) \in {\bb C}{\<}x_1,\ldots,x_n{\>}
\leqno (3.4.1)$$
of two elements given in the normal form, back into the normal form, i.e.,
expressing it as
$$\displaystyle \sum_{\nu\in\Lambda}
\ldb h_{\nu}(x)\rdb { M}_{\nu}(x).\leqno (3.4.2)$$
More precisely, we have the following fact.

\vskip .1cm

\proclaim (3.4.3) Proposition. For every
$\lambda,\mu,\nu\in \Lambda$ there is unique bilinear differential
operator with polynomial coefficients
$${ C}^{\nu}_{\lambda\mu} = C^{\nu}_{\lambda\mu} (f,g)\colon\ {\bb
C}[x_1,\ldots,x_n]\mathop{{}\otimes{}}_{\bb C} {\bb C}[x_1,\ldots,x_n]
\rightarrow {\bb
C} [x_1,\ldots,x_n]$$
with the following properties:\hfill\break
{(a)} In the normal form (3.4.2) for $\ldb f(x)\rdb
{ M}_{\lambda}(x)\cdot \ldb g(x)\rdb\cdot { M}_{\mu}(x)$ we have
$h_{\nu}
=C^{\nu}_{\lambda\mu}(f,g)$. \hfill\break
{(b)} $C^{\nu}_{\lambda\mu}=0$ if $\hbox{ord}(\nu) <
\hbox{ord}(\lambda) + \hbox{ord}(\mu)$ as well as if $\hbox{ord}(\nu)=
\hbox{ord}(\lambda)+\hbox{ord}(\mu)$ and $\nu\not=\lambda+\mu$, while
$C^{\lambda+\mu}_{\lambda\mu}$ is the operator of multiplication  of
polynomials.\hfill\break
{(c)} The $C^{\nu}_{\lambda\mu}$ satisfy the associativity constraint:
for any $
\lambda_1,
\lambda_2,
\lambda_3, \nu \in \Lambda$ the trilinear differential operators
$\displaystyle \sum_{\mu_1}C^{\nu}_{\mu_1,\lambda_3} \circ
(C^{\mu_1}_{\lambda_1,\lambda_2}\otimes 1)$ and $\displaystyle
\sum_{\mu_2}C^{\nu}_{\lambda_1,\mu_2} \circ (1\otimes
C^{\mu_2}_{\lambda_2,\lambda_3})$ coincide.

\noindent {\sl Proof:}  We first establish the existence of the
$C^{\nu}_{\lambda\mu}$ satisfying (a). In order to bring $\ldb f(x)\rdb
{ M}_{\lambda}(x) \ldb g(x)\rdb\cdot {M}_{\mu}(x)$ into the normal form,
we need
only to do this for the intermediate product, which is $M_{\lambda}(x)\ldb
g(x)\rdb$ for $\lambda\not= 0$ and $\ldb f(x)\rdb \cdot \ldb g(x)\rdb$ if
$\lambda=0$ (after this, the transforming of ${ M}_{\nu}(x){
M}_{\mu}(x)$ to
a linear combination of the ${ M}_{\sigma}(x)$ is a purely
Lie-algebraic procedure not affecting the coefficients of the form
$\ldb h(x)\rdb$).

Now, if $\lambda \not= 0$, we apply (3.3.3) to $a=M_{\lambda}(x)$ to express
$M_{\lambda}(x)\ldb g(x)\rdb$ as a sum of terms in which ordered symbols
are on the left and commutators are on the right and then use the Jacobi
identity to express each commutator via basic Lie monomials from $B$.

If $\lambda =0$, we have
$$\ldb f(x)\rdb \ldb g(x)\rdb = \ldb (f \otimes
g)({\mathop{x}^1}_{1},\ldots,
{\mathop{x}^{n}}_{n}, {\mathop{x}^{n+1}}_{1}, \ldots,
{\mathop{x}^{2n}}_n)\rdb,$$
where
$$(f\otimes g)(y_1,\ldots, y_{2n})=f(y_1,\ldots,
y_n)g(y_{n+1},\ldots,y_{2n}).$$
By using (3.3.2), we move the second copies of the $x_i$ under $f\otimes g$
to the left with the aim being to bring the second copy of $x_i$ in
adjacency with the first and use the cancellation rule (C1) of (3.1.2).
Each step of this process creates several new terms involving commutators
which, in their turn are not in the normal form, so we apply (3.3.2) to
them and so on.

This proves the existence of $C^{\nu}_{\lambda\mu}$. The uniqueness follows
since two bilinear differential operators with polynomial coefficients
whose values on every pair of polynomials coincide, should be equal as
formal expressions. Thus (a) is proved. Part (b) is clear from the nature
of the formulas (3.3.2-3).
Finally, (c) follows from the associativity of ${\bb C}{\<}x_1,\ldots,x_n{\>}$
and the fact that a trilinear differential operator with polynomial
coefficients is uniquely defined by its values on all triples of
polynomials.

\vskip .2cm

\noindent
{\bf (3.4.4) Example.} Working modulo $F^2$, i.e., setting
$$[x_i,x_j]{\varphi}[x_kx_l] = [x_i[x_j,x_k]]=0,$$
we find
$$\ldb f(x) \rdb\cdot \ldb g(x)\rdb = \ldb(fg)(x)\rdb +
\sum_{j>i}\LDB\left({\partial f\over \partial x_j}\cdot {\partial g\over
\partial x_i}\right)(x)\RDB \cdot [x_i,x_j]$$
which is equivalent to the formula of Example 1.3.9.

\vskip .1cm

\proclaim (3.4.5) Corollary. The NC-filtration on ${\bb C}
{\<}x_1,\ldots,x_n{\>}$ is given by $F^d{\bb C} {\<}x_1,\ldots, x_n{\>} =
\left\{\sum_{{\rm ord}(\lambda)\ge d} \LDB f_{\lambda} (x) \RDB
{ M}_{\lambda}(x)\right\}.$

\noindent {\sl Proof:}  Denote the LHS of the proposed equality by
$F^d$ and the RHS by $J^d$. By definition of $M_{\lambda}(x)$, we have
$J^d\subset F^d$. On the other hand, (3.4.3)(b) implies
$$J^d\cdot J^{d^\prime}\subset J^{d+d^{\prime}},\quad
[J^d,J^{d^{\prime}}]\subset J^{d+d^{\prime}+1},$$
which entails $F^d\subset J^d$.

\vskip .1cm

We now reformulate this in a more invariant form. Let $V$ be a
finite-dimentinal vector space. Consider the graded vector space
${\cal L}ie_{+} (V) = \oplus_{d\ge 2} {\cal L}ie^d(V)$ with
the grading
given by $\hbox{deg}({\cal L}ie^d(V))=d-1$. Introduce in the
symmetric algebra
$S({\cal L}ie_+ (V))$ the induced grading  and set
$$Q^d(V)=S({\cal L}ie_+(V))^d$$
to be the $d$th homogeneous part. Clearly $Q^d$ is a polynomial  functor
on the category of vector spaces [Mac].

\vskip .1cm

\proclaim (3.4.6) Proposition.  We have a natural
identification
$$Q^d(V)\simeq F^d {\cal A}ss (V)\cap V^{\otimes d}.$$
If $V=\mathop{\bigoplus}^n_{i=1}{\bb C}\cdot x_i$, so that $\hbox{$\cal
A$ss}(V)={\bb
C}{\<}x_1,\ldots, x_n{\>}$, then the ${ M}_{\lambda}(x)$,
$\hbox{ord}(\lambda)=d$,
form a basis of $Q^d(V)$.

\noindent {\sl Proof:} By construction of $F^d$, we have
$F^{d+1}{\cal A}ss (V)\cap V^{\otimes d}=0$, thus
$$F^d {\cal A}ss (V) \cap V^{\otimes d} \simeq \hbox{gr}_F^d {\cal A}ss
(V)\cap
V^{\otimes d}$$
is identified with a part of the Poisson algebra
$\hbox{gr}^{\bullet}_F{\cal A}ss (V)$. The (commutative)
multiplication and
the Poisson bracket in this algebra define a map
${\varphi}\colon\ Q^d(V)  \rightarrow \hbox{gr}^d_F {\cal A}ss(V)$
whose image is contained in the image of $F^d{\cal A}ss(V)\cap
V^{\otimes d}$, so we get
$$\psi\colon\ Q^d(V)  \rightarrow F^d {\cal A}ss (V)\cap V^{\otimes
d}.$$
To see that $\psi$ is an isomorphism, choose a basis $x_i,\ldots, x_n$ in
${V}$. Then the ${M}_{\lambda}(x), \hbox{ord}(\lambda)=d$ can be
seen as the images, under $\psi$, of the elements of the basis of
$Q^d(V)=S({\cal L}ie_+(V))^d$ formed by products of the basic Lie
monomials
in the $x_i$. On the other hand, an element
$\sum_{{\rm ord}(\lambda)\ge d}\ldb f_{\lambda}(s)\rdb M_{\lambda} (x)$
lies in $V^{\otimes d}$ if and only if $f_{\lambda}=0$ for
$\hbox{ord}(\lambda)>d$ and $f_{\lambda}\in {\bb C}$ for
$\hbox{ord}(\lambda)=d$. This proves that the $M_{\lambda}(x),
\hbox{ord}(\lambda)=d$, form a ${\bb C}$-basis in $F^d {\cal A}ss(V)\cap V^{\otimes d}$ and
that $\psi$ is an isomorphism.

\vskip .1cm

\proclaim (3.4.7) Proposition. We have a natural
identification of $GL(V)$-equivariant $S(V)$-modules:
$${\rm gr}^d_F{\cal A}ss(V)\simeq S(V)\mathop{{}\otimes{}}_{\bb C}
Q^d(V).$$

\noindent {\sl Proof:}  The identification of (3.4.6) gives a
morphism of $S(V)$-modules
$$S(V)\otimes Q^d(V)  \rightarrow \hbox{gr}^d_F {\cal A}ss(V).$$
To see that it is an isomorphism, we choose a basis in $V$ and use (3.4.5).

\vskip .1cm

Another corollary of (3.4.5) is as follows.

\proclaim (3.4.8) Proposition. The algebra ${\bb C}
{\<}x_i,\ldots
, x_n{\>}_{\ldb ab\rdb}$ can be identified with the set of possibly infinite
formal sums $\sum_{\lambda\in \Lambda}\ldb f_{\lambda} (x)
\rdb M_{\lambda}(x)$ and multiplication given by the operators $C^{\nu}_
{\lambda\mu}$.

\vskip .2cm

\noindent
{\bf (3.5) Explicit description of ${\cal O}_{A^n}^{\rm
NC}$.} Let ${\cal D}_n$ be the Weyl algebra of differential
operators in ${\bb
C}[x_1,\ldots,x_n]$. By a ${\cal D}_n$-algebra we  will mean a left ${\cal
D}_n$-module $A$ equipped with a commutative algebra structure such that
the action of the $\partial/\partial x_i$ is by algebra derivations. For
example, any localization ${\bb C}[x_i,\ldots, x_n][S^{-1}]$ is a ${\cal
D}_n$-algebra. Note that multilinear differential operators with
polynomial coefficients can be evaluated on elements of $A$. 

\vskip .1cm

\proclaim (3.5.1) Definition. Let $A$ be a ${\cal
D}_n$-algebra. The algebra $\displaystyle{\bb C} {\<}x_i,\ldots,x_n{\>}_{\ldb
ab\rdb}\otimes_{{\bb C}[x_1\ldots x_n]} A$ is defined as the space of possibly
infinite formal expressions $\sum_{\lambda\in \Lambda}\ldb
f_{\lambda}\rdb M_{\lambda}(x), f_{\lambda}\in {A}$ and the
multiplication given by
$$\ldb f_{\lambda}\rdb M_{\lambda}(x)\cdot \ldb f_{\mu}\rdb M_{\mu}
(x)=\sum_{\nu}\ldb C^{\nu}_{\lambda\mu} (f_{\lambda},f_{\mu})\rdb
M_{\mu}(x),$$
where $C^{\nu}_{\lambda\mu}$ was introduced in (3.4.3).

\vskip .1cm

\proclaim  (3.5.2) Proposition.  The above multiplication is
well defined and makes $$\displaystyle{\bb C}{\<}x_1,\ldots, x_n{\>}_{\ldb ab
\rdb}\mathop{{}\otimes{}}_{{\bb C}[x_1\ldots x_n]}A$$
 into an NC-complete
associative algebra,
whose $d$th layer of the NC-filtration consists of elements of the form
$\sum_{{\rm ord}(\lambda)\ge d} \ldb f_{\lambda} \rdb
M_{\lambda}(x)$.

\noindent {\sl Proof:}  Well-definedness: we need only to show that
the product of two infinite sums $(\sum_{\lambda}\ldb
f_{\lambda}\rdb M_{\lambda}(x))(\sum_{\mu}\ldb g_{\mu}\rdb M_{\mu}(x))$
makes sense, i.e. the number of pairs $\lambda,\mu$ such that
$C^{\nu}_{\lambda\mu}\not=0$, is finite for every $\nu$. But this follows
from (3.4.3)(b). Associativity follows from (3.4.3)(c). The statement about
the NC-filtration is proved in the same way as (3.4.5). The NC-completeness
follows from that.

\vskip .1cm

\proclaim (3.5.3) Theorem.  Let $U\subset A^n$ be a Zariski
open set, ${\cal O}(U)$ be the algebra of rational functions regular in
$\rm U$ and
$${\cal A}(U) = {\bb C}{\<}x_1,\ldots, x_n{\>}_{\ldb
ab\rdb}\mathop{{}\otimes{}}_{{\bb C}[x_1,\ldots,x_n]}{\cal O}({\rm U}).$$
Then $\Gamma(U, {\cal O}_{A^n}^{\rm NC})$ is naturally identified with
${\cal A}(U)$.

\noindent {\sl Proof:}  The ${\cal A}(U)$ obviously form a sheaf
${\cal A}$ of NC-complete algebras on $A^n$. For $d\ge 0$ let ${\cal
A}^{(d)}={\cal A}/F^{d+1}$ be the sheaf formed by expressions
$\sum_{{\rm ord}(\lambda)\le d}\ldb f_{\lambda}\rdb M_{\lambda}
(x)$ (with all the $\ldb h_{\nu}\rdb M_{\nu}(x), \hbox{ord}(\nu)>d$ in the
product of two such expression being discarded). Denote by $q_d\colon\
{\cal A}^{(d)}  \rightarrow {\cal A}^{(d-1)}$ the natural projection.
Let also $R^{(d)} ={\bb C}{\<}x_1,\ldots,x_n{\>}/F^{d+1}$ 
and ${\cal O}^{(d)}$
be
the sheaf of algebras on $A^n$ corresponding to $\hbox{Spec}(R^{(d)}) =
(A^n,{\cal O}^{(d)})$. Denote $p_d\colon\ {\cal O}^{(d)}  \rightarrow {\cal
O}^{(d-1)}$ the natural projection. To prove our theorem, it is enough to
construct, for each $d$, an isomorphism $\varphi_d\colon\ {\cal O}^{(d)}
\rightarrow {\cal A}^{(d)}$ of sheaves of algebras, in a way compatible
with the $p_d,q_d$.

\vskip .1cm

{\bf\it Construction of $\varphi_c$.} Let $U=\{f \not= 0\}$
be a principal open subset in $A^n$, so ${\cal O}(U) ={\bb C}[x_1, \ldots,
x_n][f^{-1}]$. Let $\ldb f\rdb \in {\bb C}{\<}x_1, \ldots, x_n{\>}$ be the
ordered lifting of $f$. Then $\ldb f\rdb_{ab} = f$. Note that $\Gamma (U,
{\cal O}^{(d)}) =  R^{(d)}\left[\ldb f\rdb^{-1}\right]$ (any $g$
with $g_{ab} =f$ will become invertible after $\ldb f\rdb$ is inverted). On
the other hand, we have an NC-nilpotent algebra ${\cal A}^{(d)}(U)$ with
${\cal A}^{(d)}(U)_{ab}={\cal O}(U)$, containing $R^{(d)}$. In this
algebra $\ldb f\rdb$ is invertible (because $f = \ldb f\rdb_{ab}$ is
invertible in ${\cal O}(U)$, see Proposition 2.1.1). Thus, by the universal
property of the ring of fractions, we have a homomorphism
$\varphi_{d,U}\colon\ {\cal O}^{(d)}(U)  \rightarrow {\cal A}^{(d)}(U)$.
From these we construct a morphism of sheaves $\varphi_d\colon\ {\cal
O}^{(d)}  \rightarrow {\cal A}^{(d)}$ in a standard way.
It is clear that $q_d\varphi_d=\varphi_{d-1}p_d$.

\vskip .1cm

{\bf\it $\varphi_d$ is an isomorphism.} Since $\varphi_d$
(as any homomorphism of sheaves of algebras) takes $F^i{\cal O}^{(d)}
\rightarrow F^{i}{\cal A}^{(d)}$, it suffices to prove that
$\overline{\varphi}_d\colon\ \hbox{gr}^{\bullet}_F{\cal O}^{(d)} \rightarrow
\hbox{gr}^{\bullet}_F {\cal A}^{(d)}$ is an isomorphism. As before, it is
enough
to do this over a principal open subset $U=\{f\not=0\}$. In this case, by
Theorem 2.1.6,
$$\hbox{gr}^{\bullet}_F {\cal O}^{(d)}(U)=\hbox{gr}^{\bullet}_F
\left(R^{(d)}\left[
\ldb f\rdb^{-1}\right]\right) =\hbox{gr}^{\bullet}_F (R^{(d)})[f^{-1}].$$
But $\hbox{gr}^i_F R^{(d)}$ is equal to 0 for $i>d$ and to $\hbox{gr}^i_F
{\bb C}{\<}x_1,\ldots, x_n{\>}$ for $i\le d$. The latter is, by (3.4.6), a
free
${\bb C}[x_1,\ldots, x_n]$-module with basis  $M_{\lambda}(x),
\hbox{ord}(\lambda)=i$. Therefore $\hbox{gr}^{\bullet}_F{\cal O}^{(d)}(U)$
is a free ${\cal O}(U)$-module with basis
$M_{\lambda}(x),\hbox{ord}(\lambda)\le d$. From Proposition 3.5.1,we see
that $\hbox{gr}^{\bullet}_F{\cal A}^{(d)}(U)$ has the same form. Moreover,
the morphism $\overline{\varphi}_d$ is the identity. Theorem is proved.

\vskip .2cm

\noindent {\bf (3.5.4) Examples.} (a) Let $n=2$, and ${\rm
U}\subset A^2$ be given by $x_1\not=0$. We have the (non-commuting)
elements $x_i=\ldb x_i\rdb\in {\cal O}^{\rm NC}(U)$ with $x_1$ invertible.
The two ordered quotients of $x_2$ by $x_1$ are written explicitly as
follows:
$$\displaylines{
x^{-1}_1x_2=\ldb x_2/x_1\rdb\cr
x_2x_1^{-1}=\ldb x_2/x_1\rdb +\ldb  1/x_1^2\rdb [x_1,x_2]+\ldb 1/x^3_1\rdb
[x_1,[x_1,x_2]]+\cdots\cr}$$
See (2.1.7.2) and (3.3.3).

\vskip .1cm

(b) Let $n=2$ and $U\subset A^2$ given by $x_1+x_2\not=0$. Identifying, as
before, $x_i$ with $\ldb x_i\rdb \in {\cal O}(U)$, we find that $x_1 + x_2$
is invertible in ${\cal O}(U)$. The inverse can be written explicitly by
applying the Taylor formula (3.3.4) to $f(z)=1/z$, getting
$$(x_1 + x_2)^{-1} = \sum^{\infty}_{k=0}\Bigldb {(-1)^kk!\over
\displaystyle {\mathop{x}^1}_1+{\mathop{x}^3}_2}\cdot
{\mathop{X}^2}_k\Bigrdb$$
Then we should bring each summand to the normal form by applying (3.3.2) to
symbols like
$$g(y_1,y_2, y_3)={y_2\over (y_1+y_3)^k}.$$

(c) Let $n=m^2$, so $A^n=\hbox{Mat}_m({\bb C})$ is the space of matrices.
The matrix elements $x_{ij}$ are the coordinates in $A^n$. Let ${\rm U} =
GL_m({\bb C})\subset A^n$ be given by $\hbox{det}\| x_{ij}\|\not=0$.
Let us now identify the $x_{ij}$ with the noncommuting elements
$\ldb x_{ij}\rdb \in {\cal O}(U)$, By (2.1.1), the tautological matrix
${\bf M}=\| x_{ij}\|\in \hbox{Mat}_m ({\cal O}(U))$ is invertible. So the entries
$C_{ij}$ of ${\bf M}^{-1}$ are some series $\sum\ldb f_{\lambda}(x)\rdb
M_{\lambda}(x), f_{\lambda} \in {\cal O} (GL_m ({\bb C}))$.

\vskip .2cm

\proclaim (3.5.5) Theorem. The algebra ${\bb
C}{\<\!\<}x_i,\ldots, x_n{\>\!\>}$ is naturally identified with $$\displaystyle{\bb
C}{\<}x_i,\ldots,x_n{\>}_{\ldb ab\rdb}\otimes_{{\bb C}[x_1\ldots x_n]}{\bb
C}[[
x_1\ldots x_n]]. $$

{\bf\it Proof.}\enspace\enspace Any formal sum $\sum_{\lambda\in \Lambda}\ldb
f_{\lambda}(x)\rdb M_{\lambda}(x)$ with $f_{\lambda}(x) \in {\bb
C}[[ x_i,\ldots, x_n]]$, can be regarded as a noncommutative power series
in the $x_i$, so we have an embedding from one algebra to another. This
embedding induces an isomorphism on  $\hbox{gr}^{\bullet}_{F}$. Both
algebras being NC-complete, our embedding is an isomorphism.

\vfill\eject

\centerline {\bf \S4. Detailed study of algebraic NC-manifolds.}

\vskip 1cm

\noindent
{\bf (4.1) Poisson envelopes.} Let ${\cal P}ois$ be the category of Poisson
algebras (1.1) and $c: {\cal P}ois \to {\cal C}om$ the functor forgetting
the
Poisson bracket.

\proclaim (4.1.1) Proposition. The functor $c$ has a left adjoint $P$
called the Poisson envelope.

\noindent {\sl Proof:} 
Let $A$ be a commutative algebra and ${\cal P}ois(A)$ be the
free Poisson algebra generated by $A$ as a vector space.  Let $\cdot$,
$\{$, $\}$ be the operations in ${\cal P}ois(A)$.  Define $P(A)$ as the
quotient of ${\cal P}ois(A)$ by the Poisson ideal generated by $a \cdot b -
ab$, $a, b \in A$.  Here $ab \in A$ is the product of $a$ and $b$ in
$A$.  The adjointness of $P$ and $c$ is obvious from the construction.

\vskip .2cm

\noindent {\bf (4.1.2) Example.} If $A = {\bb C}[x_1, \ldots, x_n]$ is the
free commutative algebra, then $P(A) = {\cal P}ois (x_1, \ldots, x_n)$ is the
free Poisson algebra on $x_1, \dots, x_n$. By (3.4.5-6) this algebra is
canonically identified with ${\rm gr}^{\bullet}_F {\bb C}\langle x_1, \ldots,
x_n\rangle$.

\vskip .2cm

Note that by construction $P(A) = \mathop{\bigoplus}\limits_{d \geq 0}
P^d(A)$ is a graded (1.1.3) Poisson algebra, with $P^d(A)$ spanned by
the images of formal expressions from $\hbox{Pois}(A)$ containing $d$
instances of Poisson brackets.  In particular $P^0(A)=A$ and each
$P^d(A)$ is an $A$-module.

In (3.4) we defined a polynomial functor $Q^d$ on the category of
vector spaces.  Let $Q_A^d$ be the natural extension of $Q^d$ to the
category of projective $A$-modules.  For such a module $M$ the fiber
of $Q_A^d(M)$ at a ${\bb C}$-point $x \in \hbox{Spec}(A)$ is the value of
$Q^d$ on the fiber of $M$ at $x$.

\proclaim (4.1.3) Theorem. Let $A$ be a smooth, finitely generated,
commutative algebra.  Then $P^d(A) \simeq Q_A^d(\Omega_A^1)$ as an
$A$-module.

\noindent {\sl Proof:} We first construct an $A$-linear map $\varphi: P^d(A)
\to Q_A^d(\Omega_A^1)$. The space $P^d(A)$ is spanned by products of
the form $\pi = f_0 h_1 \cdots h_m$, where $f_0 \in A$ and each $h_v$
is a bracket monomial
$$h_v = \{f_{v,0}, \{f_{v,1}, \ldots, \{f_{v,l_v-1},
f_{v,l_v}\}\ldots\}, \qquad f_{v,j} \in A, \qquad \sum_{v=1}^m l_v =
d.$$
These products are considered modulo the Jacobi and Leibniz identities.
For $\pi$ as above set
$$\varphi(\pi) = f_0 \cdot \prod_{v=1}^m [df_{v,0}, [df_{v,1}, \ldots,
[df_{v,l_v-1}, df_{v,l_v}] \ldots ] \in Q_A^d(\Omega_A^1).$$
Because $d: A \to \Omega_A^1$ satisfies the Leibniz rule, this is
compatible with the Leibniz identity for the Poisson brackets.
Because [ , ] satisfies the Jacobi identity, $\varphi$ is compatible
with the Jacobi identity for Poisson brackets.  So $\varphi$ is a well
defined $A$-module homomorphism.

We now prove that $\varphi$ is an isomorphism.  Both the Poisson
envelope and localization being direct limit-type constructions, we
have
$$P^d(A[S^{-1}]) = P^d(A)[S^{-1}]. \leqno{(4.1.4)}$$
Thus it is enough to prove that $\varphi$ is an isomorphism in the
case when $A$ is a local ring, obtained by localizing a finitely
generated smooth algebra at a ${\bb C}$-point.  In this case we
represent
$$A = {\bb C}[x_1, \ldots, x_n]^{\sim}/(h_1, \ldots, h_m)$$
where the superscript 
``$\sim$'' means localization at 0 and $h_1, \ldots, h_m \in {\bb
C}[x_1, \ldots, x_n]$ are such that $d_0h_i$ are linearly independent.
Then, denoting by $\{(h_1, \ldots, h_m)\} \subset {\cal P}ois (x_1, \ldots,
x_n)$ the Poisson ideal generated by the $x_i$, we have, in virtue of
(4.1.4):
$$P^d(A) \cong ({{\cal P}ois}(x_1, \ldots, x_n)/\{(h_1, \ldots, h_m)\})^d
\mathop{\otimes}\limits_{{\bb C}[x_1 \ldots x_n]} {\bb C}[x_1, \ldots,
x_n]^{\sim}. \leqno{(4.1.5)}$$
Now, since
$$
{\cal P}ois^d(x_1, \ldots, x_n) \simeq {\bb C}[x_1, \ldots, x_n]
\mathop{\otimes}\limits_{{\bb C}} Q^d(\oplus {\bb C} \cdot x_i) \simeq$$
$$\simeq  Q^d_{{\bb C}[x_1, \ldots, x_n]} (\Omega^1_{{\bb C}[x_1, \ldots,
x_n]})
\leqno (4.1.6)$$
and
$$\Omega_A^1 \simeq \Omega^1_{{\bb C}[x_1, \ldots, x_n]^{\sim}} / (h_i,
dh_i) \leqno{(4.1.7)}$$
we find that the factorization of ${{\cal P}ois}^d(x_1, \ldots, x_n)$ given by
$\{(h_1, \ldots, h_m)\}$ exactly corresponds, via $\varphi$, to the
factorization of $Q^d_{{\bb C}[x_1, \ldots, x_n]}(\Omega^1_{{\bb
C}[x_1, \ldots, x_n]})$ induced by (4.1.7).  So $\varphi$ is an
isomorphism and the theorem is proved.

\vskip .1cm

Let now $M$ be a smooth algebraic variety.  Because of (4.1.4), we
have a sheaf $P({\cal O}_M)$ of graded Poisson algebras on $M$ with
$P^0({\cal
O}_M)
= {\cal O}_M$ and $P^d({\cal O}_M) \simeq Q^d(\Omega_M^1)$ as an ${\cal
O}_M$-module.  The
definition of $Q^d$ implies:

\proclaim  (4.1.8) Proposition. The spectrum of $P({\cal O}_M)$ is the
(infinite-dimensional) manifold ${\mathop\prod\limits_{d\geq2}}{}_{_M}
{{\cal L}ie}^d(T_M)$, i.e. the completion of the degree${}\geq 2$ part of the
free ${\cal O}_M$-Lie algebra generated by the tangent bundle $T_M$.

\vskip .1cm

\noindent {\bf (4.1.9) Example.} A Poisson structure $\xi$ on $M$ itself
gives a section $M \to \hbox{Spec}(P({\cal O}_M))$.  Its composition with
the
projection to ${{\cal L}ie}^2(T_M) = \bigwedge ^2 T_M$ is the bivector
corresponding to $\xi$ in a standard way.

\vskip .3cm

\noindent {\bf (4.2) Structure of ${\rm gr}^{\bullet}_F$ for $d$-smooth algebras.}

\proclaim (4.2.1) Theorem. Let $R$ be a $d$-smooth algebra, $A =
R_{ab}$.  Then the Poisson algebra ${\rm gr}^{\bullet}_F(R)$ is canonically
isomorphic to $P(A)/P^{\geq d+1}(A)$.  In particular, ${\rm gr}^{i}_F(R)
\simeq Q_A^i(\Omega_A^1)$ for $i \leq d$.

\noindent {\sl Proof:} Since ${\rm gr}^{\bullet}_F(R)$ is a graded Poisson algebra
containing ${\rm gr}^0_F(R) = A$, the universal property of $P(A)$ gives a
morphism $P(A) \to {\rm gr}^{\bullet}_F(R)$ vanishing on $P^{\geq d+1}(A)$, so
we get $\varphi: P(A)/P^{\geq d+1}(A) \to {\rm gr}^{\bullet}_F(R)$.  We will
prove that $\varphi$ is an isomorphism.  In fact, by induction we can
assume that the homogeneous part $\varphi^i: P^i(A) \to {\rm gr}^i_F(R)$ is
an isomorphism for $i < d$, so only $\varphi^d$ needs to be treated.
Note that ${\rm gr}^d_F(R) = F^dR$ is a central ideal in $R$.

Let $x \in \hbox{Spec}(A)$ be a ${\bb C}$-point.  For an $A$-module $N$ we
denote $N|_x = N/{\bf m}_{x,ab}N$ the fiber of $N$ at $x$.  Here ${\bf m}_{x,ab}
\subset A$ is the ideal of $x$.  Since the source and the target of
$\varphi^d$ are $A$-modules of finite type, it is enough to prove that
for any $x$ as above the map of fibers $\varphi^d|_x: P^d(A)|_x \to
(F^dR)|_x$ is an isomorphism.  Let $(R_x, {\bf m})$ and $(A_x, {\bf m}_{ab})$ be the
localizations of $R$ and $A$ at $x$.  Since $P^d$ and ${\rm gr}^d_F$ commute
with localizations, we can identify $P^d(A)|_x \cong P^d(A_x)|_x$,
$(F^dR)|_x \cong (F^dR_x)|_x$.  Consider the commutative diagram of
$A_x$-modules
$$\matrix{%
P^d(A_x) & \mathop{\to}\limits^{\varphi^d} & F^d(R_x)\cr
\llap{$\scriptstyle\alpha$} \downarrow & & \downarrow
\rlap{$\scriptstyle\beta$}\cr
P^d(A_x/{\bf m}_{ab}^j) & \mathop{\to}\limits^{\psi} & F^d(R_x/{\bf m}^j)}$$
where $j\gg0$, the maps $\alpha, \beta$ are induced by the
functoriality of $P^d$ and $F^d$, and $\psi$ is defined similarly to
$\varphi$, via the universal property of Poisson envelopes.  We claim
that for $j\gg0$ the maps $\alpha,\beta,\psi$ give isomorphisms on
fibers at $x$.

The fact that $\alpha$ gives an isomorphism follows from (4.1.3) and
from an identification $A_x/{\bf m}_{ab}^j \simeq {\bb C}[x_1, \ldots, x_n] /
(x_1, \ldots, x_n)^j$.  Namely, we explicitly identify the Poisson
envelope of the latter algebra as a quotient of $P({\bb C}[x_1, \ldots,
x_n]) = {\cal P}ois(x_1, \ldots, x_n)$.

The fact that $\psi$ gives an isomorphism, follows from the
identification
$$R_x/{\bf m}^j \simeq {\bb C}\langle x_1, \ldots, x_n\rangle \biggl/(x_1, \ldots,
x_n)^j + F^{d+1},$$
see (1.5.1), and from the explicit description of $F^{\bullet}{\bb
C}\langle x_1, \ldots, x_n\rangle$ given in (3.4.5).

Finally, we prove that $\beta$ gives an isomorphism on fibers.  We
have
$$\eqalign{%
F^d(R_x)|_x &= F^d(R_x)/{\bf m} \cdot F^d(R_x)\cr
F^d(R_x/{\bf m}^j)|_x &= F^d(R_x)/{\bf m} \cdot F^d(R_x) + {\bf m}^j \cap F^d(R_x)\cr}$$
So we are reduced to the next lemma.

\proclaim (4.2.2) Lemma. If $j > d$, then ${\bf m}^j \cap F^d(R_x) \subset
{\bf  m}
\cdot F^d(R_x)$.

\noindent {\sl Proof:} Let $(\hat{R}_x, \hat{\bf m})$ be the $m$-adic completion of
$R$.  By (1.5.1),
$$\hat{R}_x \simeq {\bb C}\langle\!\langle x_1, \ldots,
x_n\rangle\!\rangle/F^{d+1}.$$
This is again a Noetherian local algebra.  Now, the explicit
description of ${\bb C}\langle\!\langle x_1, \ldots, x_n\rangle\!\rangle$
together with its NC-filtration given in (3.5.5), implies that
$\hat{\bf m}^j \cap F^d\hat{R}_x \subset \hat{\bf m} \cdot F^d\hat{R}_x$ for
$j > d$.  Our statement would follow therefore from the equality
$${\bf m} \cdot F^d(R_x) = (\hat{\bf m} \cdot F^d(\hat{R}_x)) \cap R
\leqno{(4.2.3)}$$
But this is a particular case of the following version of the Krull
lemma whose proof is the same as in the classical case ([Eis], Ch. 4).

\proclaim (4.2.4) Lemma. Let $R$ be a (not necessarily commutative)
Noetherian local ring, $\hat{R}$ is completion, $I \subset R$ a
central ideal, $\hat{I} = I\hat{R}$.  If $\hat{I}$ is central, then
$\hat{I} \cap R = I$.

Theorem 4.2.1 is proved.

\vskip .1cm

\proclaim (4.2.5) Corollary. Let $R$ be a $d$-smooth algebra.  Then
$H_2(R,R_{ab}) \simeq Q^{d+1}_{R_{ab}}(\Omega^1_{R_{ab}})$ is an
$R_{ab}$-module.

\vskip .1cm

\noindent {\bf (4.2.6) Remark.} It would be interesting to describe
$H_i(R,R_{ab})$, $i > 2$ as certain explicit polynomial functors of
$\Omega^1_{R_{ab}}$, generalizing the Hochschild-Kostant-Rosenberg
theorem $H_i(R_{ab},R_{ab}) = \Lambda^i \Omega^1_{R_{ab}}$.  The
analogous group-theoretical problem (find the integral homology of the
free $(d+1)$-stage nilpotent group, see (1.6.4)) is quite difficult.  Even rationally, when groups can be replaced by Lie algebras,
the answer for the full homology has been found only for $d \leq 1$,
see [Sig].

\vskip .3cm

\noindent {\bf (4.3) Cohomological classification of thickenings.} Let ${M}$
be
a smooth algebraic variety and $X \supset { M}$ be a $d$-smooth
thickening.  For a Zariski open ${ U} \subset {\cal M}$ let ${
U}^{(d)}$ be the induced thickening of $U$.  Let ${\cal T}_X(U)$ be
the category whose objects are $(d+1)$-smooth thickenings $V \supset
U^{(d)}$ and morphisms are morphisms of NC-schemes identical on
$U^{(d)}$.  Denote by $T_M$ the tangent sheaf of $M$.

\vskip .1cm

\proclaim  (4.3.1) Proposition. (a) Each ${\cal T}_X(U)$ is a groupoid.\hfill\break
(b) The correspondence $U \mapsto {\cal T}_X(U)$ is a stack ${\cal T}_X$ of
groupoids on the Zariski topology of ${\cal M}$.\hfill\break
(c) The stack ${\cal T}_X$ is a gerbe with band
$\hbox{Hom}(Q^{d+1}T_M,T_M)$.

\noindent {\sl Proof:} (a) follows from (1.2.7).  Part (b) is clear.  Part (c)
follows from (4.2.5) and (1.3.7).

\vskip .1cm

By applying the Grothendieck-Giraud formalism of nonabelian cohomology via
stacks [Bry], we get the following result.

\proclaim (4.3.2) Theorem. Any $d$-smooth thickening $X \supset { M}$
defines a natural cohomology class (obstruction) $\gamma_X \in H^2(M,
\hbox{Hom}(Q^{d+1}T_M,T_M))$.  A $(d+1)$-smooth thickening ${Y}
\supset X$ exists if and only if $\gamma_X = 0$.  If this is the case,
the set of such $Y$ modulo isomorphisms identical on $X$, is a
principal homogeneous space over $H^1(M, \hbox{Hom}(Q^{d+1}T_M,T_M))$.

\vskip .2cm

\noindent {\bf (4.4) Jet Bundles.} Let ${ M}$ be a smooth algebraic variety,
$\dim({ M}) = n$.  For $x \in { M}$ let ${\bf m}_{x, ab} \subset {\cal
O}_{
M}$ be the ideal
of ${ M}$.  Denote also by ${\bf m}_{ab} \subset {\bb C}[x_1, \ldots, x_n]$
the
ideal $(x_1, \ldots, x_n)$, and by $\bf m$ the similar ideal in ${\bb
C}\langle x_1, \ldots, x_n\rangle$.

Let $l \geq 0$ and $J_{ M}^l$ be the bundle (sheaf) of $l$-jets of
functions on ${ M}$.  The fiber of $J_M^l$ at $x \in M$ is
${\cal O}_M/{\bf m}_{x,ab}^{l+1}$.  Thus $J_{\cal M}^l$ is a sheaf of ${\cal
O}_{\cal M}$-algebras
locally isomorphic to ${\cal O}_{\cal M} \otimes ({\bb C}[x_1, \ldots,
x_n]/{\bf m}_{ab}^{l+1})$.  Let $j_l^{ab}: {\cal O}_M \to J_{ M}^l$ be the
universal
$l$th order differential operator.  The sections of $\hbox{Im}(j_l^{ab})$
will be called pure jets.

Let now $X$ be a $d$-smooth thickening of ${\cal M}$.  For $x \in {\cal M}$
let ${\bf m}_x
\subset {\cal O}_X$ be the preimage of ${\bf m}_{x,ab}$.  We have then a bundle
(locally free sheaf) $J_X^l$ on $M$ whose fiber at $x \in X$ is
${\cal O}_X/{\bf m}_x^{l+1}$.  Again, we have a natural morphism of sheaves $j_l:
{\cal O}_X \to J_X^l$, lifting $j_l^{ab}$.

\vskip .1cm

\proclaim (4.4.1) Proposition. $J_X^l$ is a sheaf of ${\cal O}_M$-algebras
locally isomorphic to 
$${\cal O}_M \otimes ({\bb C}\langle x_1, \ldots,
x_n\rangle/({\bf m}^{l+1} + F^{d+1})).$$
  Its abelianization is identified with
$J_M^l$.

\noindent {\sl Proof:} follows from (1.5.1).

\vskip .1cm

Let ${ U} \subset M$ be Zariski open, and $U^{(d)} \subset X$ be
the induced thickening.  Define the following category ${\cal
J}_X(U)$.  Its objects are pairs $({\bb O},\psi)$, where ${\bb O}$ is
a sheaf of ${\cal O}_U$-algebras locally isomorphic to ${\cal O}_U \otimes
({\bb
C}\langle x_1, \ldots, x_n\rangle/{\bf m}^{d+2})$, and $\psi: {\bb O} \to
{\cal J}_X^{d+1}$ is a surjection of sheaves of algebras with kernel
$F^{d+1}{\bb O}$.  A morphism $({\bb O},\psi) \to ({\bb O}',\psi')$ is
an isomorphism $f: {\bb O} \to {\bb O}'$ of sheaves of algebras such
that $\psi'f = \psi$.

\noindent Note that if ${\bb O}$ is any sheaf of ${\cal O}_U$-algebras locally
isomorphic to ${\cal O}_U \otimes ({\bb C}\langle x_1, \ldots,
x_n\rangle/{\bf m}^{d+2})$, then the images of ${\cal O}_U \otimes 
({\bf m}/{\bf m}^{d+2})$
under all such isomorphisms coincide and define a sheaf of ideals
$\underline{\bf m} \subset {\bb O}$, with ${\bb O}/\underline{\bf m} = {\cal O}_U$.
Define the category ${\cal A}_X(U)$ whose objects are sheaves of
algebras ${\bb O}$ as before together with an isomorphism $\varphi:
{\bb O}/{\bf m}^{d+1} \mathop{\to}\limits^{\sim} J_X^d$ and morphisms
defined in a similar way.  We have the functors
$${\cal T}_X(U)\mathop{\to}\limits^{e_U} {\cal J}_X({\bf U})
\mathop{\to}\limits^{f_U} {\cal A}_X(U)$$
defined as follows.  For a $(d+1)$-smooth thickening ${V} \supset
U^{(d)}$ we set $e_U(V) = {\cal J}_{V}^{d+1}$.  For an object
$({\bb O},\psi)$ of ${\cal J}_X(U)$ we set $f_U({\bb O},\psi) = ({\bb
O},p\psi)$ where $p: J_X^{d+1} \to J_X^d$ is the natural surjection.

\vskip .1cm

\proclaim (4.4.2) Proposition. (a) The correspondences $U \mapsto {\cal
J}_X(U), {\cal A}_X(U)$ form stacks of groupoids ${\cal J}_X, {\cal
A}_X$ on $M$, and the functors $e_U, f_U$ define morphisms of stacks
${\cal T}_X \mathop{\to}\limits^e {\cal J}_X \mathop{\to}\limits^f {\cal
A}_X$.\hfill\break
(b) The stack ${\cal J}_X$ is a gerbe with band
$\hbox{Hom}(Q^{d+1}T_M,T_M)$,
the stack ${\cal A}_X$ is a gerbe with band
$\hbox{Hom}(T_M^{\otimes(d+1)},T_M)$.\hfill\break
(c) The functor $e$ is an equivalence of gerbes with band
$\hbox{Hom}(Q^{d+1}T_M,T_M)$.\hfill\break
(d) The functor $f$ is compatible with the morphism 
$${\epsilon}:
\hbox{Hom}(Q^{d+1}T_M,T_M) \to \hbox{Hom}(T_M^{\otimes(d+1)},T_M)$$
 induced
by the
embedding $Q^{d+1}\Omega^1_M \hookrightarrow (\Omega^1_M)^{\otimes
d+1}$.

\noindent {\sl Proof:} (a) is clear.  Let us prove (b).  Let first an object
$({\bb O},\varphi)$ of ${\cal A}_X(U)$ be given.  Then
$\underline{\bf m}/\underline{\bf m}^2 \simeq \Omega^1_U$ and thus
$\underline{\bf m}^i/\underline{\bf m}^{i+1} \simeq (\Omega^1_U)^{\otimes i}$
for $i \leq d+1$, while $\underline{\bf m}^{d+2} = 0$.  Thus ${\bb O}$ is
a central extension
$$0 \to (\Omega^1_M)^{\otimes d+1} \to {\bb O}
\mathop{\to}\limits^{\varphi} {\bb O}/\underline{\bf m}^{d+1} \to 0$$
Automorphisms of such an extension compatible with $\varphi$, form a
Abelian group sheaf of Abelian groups 
identified with $\hbox{Hom}(T_m^{\otimes d+1},T_M)$, see
(1.2.5).  So ${\cal A}_X$ is a gerbe with band
$\hbox{Hom}(T_M^{\otimes(d+1)},T_M)$.

Further, an object $({\bb O},\psi)$ of ${\cal J}_X(U)$ gives a central
extension
$$0 \to F^{d+1}{\bb O} \to {\bb O} \mathop{\to}\limits^{\psi} {\cal
J}_X^{d+1} \to 0$$
and our statement follows in the same way as above once we recall the
interpretation of the functor $Q^{d+1}$ given in (3.4.6).

To prove (c) it is enough to show that $e$ is just a morphism of
gerbes with the same band: it will then be automatically an
equivalence.  But this is clear from the definition of $e$: given
${\bb O} = J_{ U}^{d+1}$, the consideration of jets at any $x \in
M$ gives a canonical identification $F^{d+1}{\bb O}
\mathop{\leftarrow}\limits^{\sim} F^{d+1}{\bb O}_{ U}$.  As for
(d), this just expresses the fact that for $({\bb O},\psi) \in {\cal
J}_X(U)$ the embedding $\hbox{Aut}({\bb O},\psi) \hookrightarrow
\hbox{Aut}({\bb
O},p\psi)$ is precisely identified with the embedding
$$\hbox{Hom}(Q^{d+1}T_M,T_M) \to \hbox{Hom}(T_M^{\otimes(d+1)},T_M)$$
mentioned in the statement.  But this is also clear, since
$F^{d+1}{\bb O} \subset \underline{\bf m}^{d+1}$ is precisely
$Q^{d+1}(\Omega^1_M) \hookrightarrow (\Omega^1_M)^{\otimes(d+1)}$.

\vskip .2cm

Let us reformulate part (c) of the above proposition.

\proclaim  (4.4.3) Corollary. In order to extend a given $d$-smooth
thickening $X \supset M$ to a $(d+1)$-smooth thickening, it suffices
to construct a sheaf ${\bb O}$ of ${\cal O}_M$-algebras on $M$ locally
isomorphic to ${\cal O}_M \otimes ({\bb C}\langle x_1, \ldots,
x_n\rangle / {\bf m}^{d+2})$ such that ${\bb O}/F^{d+1}{\bb O} \simeq
J_X^{d+1}$.

Note that in this case ${\bb O}_{ab} \simeq J_M^{d+1}$, and the
structure sheaf ${\cal O}_{ {Y}}$ of the thickening consists of
those sections $g$ of ${\bb O}$ for which $g_{ab}$ is a pure jet.

\vskip .1cm

Let ${\cal A}ss(\Omega^1_M)$ be the tensor algebra of $\Omega^1_M$
over ${\cal O}_M$ and ${\cal A}ss^{\leq d}(\Omega^1_M)$ its quotient
obtained by disregarding the homogeneous components of degree
$>d$. Let ${\cal G}_d$ be the sheaf of groups on $M$ formed by algebra
automorphisms of ${\cal A}ss^{\leq d}(\Omega^1_M)$ whose action on
$\Omega^1_M$ is identical modulo $(\Omega^1_M)^{\otimes i}$, $i \geq
2$.  Similarly, let $S^{\leq d}(\Omega^1_M)$ be the truncated
symmetric algebra and ${\cal G}'_d$ be the sheaf of automorphisms of
$S^{\leq d}(\Omega^1_M)$ identical modulo terms of degree $\geq 2$.
Clearly we have a surjective homomorphism $\sigma_d: {\cal G}_d \to
{\cal G}'_d$.  The next proposition is straightforward.

\proclaim  (4.4.4) Proposition. (a) ${\cal G}_d$ and ${\cal G}'_d$ are
sheaves of nilpotent groups on $M$ fitting into short exact sequences
with central kernels:
$$
1 \to \hbox{Hom}(T_M^{\otimes d+1},T_M) \to {\cal G}_{d+1} \to {\cal G}_d
\to 1 \leqno (4.4.4.1)$$
$$1 \to \hbox{Hom}(S^{d+1}T_M,T_M) \to {\cal G}'_{d+1} \to {\cal G}'_d \to 1.
\leqno 
(4.4.4.2)$$
(b) $H^1(M,{\cal G}_d)$ is identified with the set of isomorphism
classes of sheaves of algebras ${\bb O}$ on $M$ locally isomorphic to
${\cal O}_M \otimes ({\bb C}\langle x_1, \ldots, x_n\rangle/{\bf m}^{d+1})$,
equipped with an identification $\zeta: \underline{\bf m}/\underline{\bf m}^2
\to \Omega^1_M$ (modulo isomorphisms preserving these
identifications).\hfill\break
(c) Similarly, $H^1(M,{\cal G}'_d)$ is identified with the set of
isomorphism classes of sheaves ${\bb O}'$ of ${\cal O}_M$-algebras locally
isomorphic to ${\cal O}_M \otimes ({\bb C}[x_1, \ldots, x_n]/{\bf m}_{ab}^{d+1})$
equipped with an identification
$\underline{\bf m}_{ab}/\underline{\bf m}^2_{ab} \simeq \Omega^1_M$.

Note, in particular, that the sheaf $J_M^d = {\bb O}'$ satisfies the
conditions of (c) so gives same class $\alpha_{d,M} \in H^1(M,{\cal
G}'_d)$.  Further, any $d$-smooth thickening $X \supset M$ gives the
sheaf $J_X^d = {\bb O}$ satisfying the conditions of (b), so gives a
class $\alpha_{d,X} \in H^1(M,{\cal G}_d)$ such that
$\sigma_d(\alpha_{d,X}) = \alpha_{d,M}$.

\vskip .1cm

\proclaim  (4.4.5) Proposition. If $X$ is a $d$-smooth thickening of
$M$, then the image of the obstruction $\gamma_X \in H^2(M,
\hbox{Hom}(Q^{d+1}T_M,T_M))$ under
$${\epsilon}_*: H^2(M, \hbox{Hom}(Q^{d+1}T_M,T_M)) \to H^2(M,
\hbox{Hom}(T_M^{\otimes(d+1)}T_M,T_M))$$
is equal to $\delta(\alpha_{d,X})$ where $\delta$ is the nonabelian
coboundary map associated to (4.4.4.1).

\noindent {\sl Proof:} This is a consequence of 4.4.2(d).

\vskip .1cm

\proclaim (4.4.6) Proposition. The map $\epsilon$ from 4.4.2(d) is the
embedding of a direct summand.  In particular the induced map $\epsilon_*$ on the cohomology is injective.

\noindent {\sl Proof:} For any $d$ the value of the polynomial functor $Q^d$ on
a vector space $V$ can be described as
$$Q^d(V) = (V^{\otimes d} \otimes W_d)^{S_d}$$
where $W_d$ is a certain representation of the symmetric group $S_d$.
The embedding $Q^d(V) \subset V^{\otimes d}$ corresponds to an
embedding $W_d \subset {\bb C}[S_d]$ of $W_d$ into the regular
representation.  Now taking the complement $W'_d$ to $W_d$ in ${\bb
C}[S_d]$ and defining $R^d(V) = (V^{\otimes d} \otimes W'_d)^{S_d}$,
we find that $V^{\otimes d} = Q^d(V) \oplus R^d(V)$ as ${GL}(V)$-modules.
Thus $(\Omega^1_M)^{\otimes d} \simeq Q^d(\Omega^1_M)
\oplus R^d(\Omega^1_M)$ as vector bundles, whence the statement.

\vskip .3cm

\noindent {\bf (4.5) First order thickenings and Atiyah classes.} Let $M$ be
as before.  Among 1-smooth thickenings of $M$ there is a distinguished
one $X_0$, with ${\cal O}_{X_0} = {\cal O}_M \oplus \Omega^2_M$ and the
multiplication given in (1.3.9).  Recalling that the functor $Q^1 =
\Lambda^2$ is the second exterior power, we get from (4.3.2) the
following.

\proclaim (4.5.1) Proposition. The set of isomorphism classes of
1-smooth thickenings of $M$ is canonically identified with
$H^1(M,\Omega^2_M \otimes T_M)$, so that $X_0$ corresponds to 0.

\vskip .1cm

We denote by $\alpha_X^- \in H^1(M,\Omega^2_M \otimes T_M)$ the class
corresponding to a thickening $X$.  Recall (4.4) that we have the
classes
$$\eqalign{%
\alpha_{2,X} \in H^1({\cal M}, {\cal G}_2) &= H^1({\cal M},
\hbox{Hom}(T_M^{\otimes
2},T_M)),\cr
\alpha_{2,M} \in H^1({\cal M}, {\cal G}'_2) &= H^1({\cal M},
\hbox{Hom}(S^2T_M,T_M))}$$
classifying the bundles ${\cal J}_X^2$, ${\cal J}_M^2$.  In fact,
$\alpha_{2,M}$ is the Atiyah class of the tangent bundle of $M$, see
[Kap].  So we call $\alpha_{2,X}$ the NC-Atiyah class of $X$.

\vskip .1cm

\proclaim (4.5.2) Proposition. With respect to the natural
identification $$H^1(M, \hbox{Hom}(T_M^{\otimes 2},T_M)) = H^1(M,
\hbox{Hom}(S^2T_M,T_M))
\oplus H^1(M,\Omega^2_M \otimes T_M)$$ we have $\alpha_{2,X} =
\alpha_{2,M} \oplus \alpha_X^-$.

\noindent {\sl Proof:} $\Omega^2_M \otimes T_M$ is the sheaf of automorphisms of
$X_0$ identical on $M$.  The identification of (4.5.1) is obtained
explicitly by identifying, for any affine ${ U} \subset M$, the
induced thickening ${ U}^{(1)} \subset X$ with the trivial
thickening ${ U}_0 \subset X_0$ and looking  at the discrepancies over
the intersections which give a \v{C}ech 1-cocycle in $\Omega^2_M
\otimes T_M$ representing $\alpha_X^-$.  Similarly, $\alpha_{2,X}$ is
obtained by identifying $J_X^2|_U$ with ${\cal A}ss^{\leq 2}
(\Omega^1_{ U})$ by an algebra automorphism with identical linear
part, and then taking the discrepancies.

Notice now that ${\cal A}ss^{\leq 2}(\Omega^1_U) = S^{\leq
2}(\Omega^1_U) \oplus \Omega^2_U$.  So if it is identified with
$J_X^2|_U$, then ${\cal O}_{{ U}^{(1)}}$ (realized inside $J_X^2|_U$ as
the sheaf of sections whose image in $J_U^2$ is a pure jet) becomes
explicitly identified with the direct sum of ${\cal O}_{ U} \subset
J_{\bf U}^2 \simeq S^{\leq 2}(\Omega^1_U)$ and $\Omega^2_U$.  In
other words, ${U}^{(1)}$ becomes identified with ${ U}_0$.
This shows that the antisymmetrization of $\alpha_{2,X}$ is
$\alpha_X^-$, proving the proposition.

\vskip .1cm

\proclaim (4.5.3) Theorem. Let $X$ be a 1-smooth thickening of $M$ and
$$\alpha = \alpha_{2,X} \in H^1(M, \hbox{Hom}(T_M^{\otimes 2},T_M))$$
be its
NC-Atiyah class.  Then $X$ extends to a second order thickening if and
only if $\alpha$ satisfies the anti-associativity condition:
$$\alpha \circ (\alpha \otimes 1) = -\alpha \circ (1 \otimes
\alpha)\quad {\rm in}\quad H^2(M, \hbox{Hom}(T_M^{\otimes 3},T_M)).$$

\noindent {\sl Proof:} The sheaf ${\cal G}_2$ is a central extension
$$1 \to \hbox{Hom}(T_M^{\otimes 3},T_M) \to {\cal G}_2 \to
\hbox{Hom}(T_M^{\otimes
2},T_M) \to 0. \leqno{(4.5.4)}$$
By (4.4.5-6) the extension of $X$ to a 2-smooth thickening is possible
if and only if $\delta(\alpha) = 0$, where $\delta$ is the coboundary
map of (4.5.4).  So we need to prove that
$$\delta(\alpha) = \alpha \circ (\alpha \otimes 1) + \alpha \circ (1
\otimes \alpha). \eqno{(4.5.5)}$$
For this, let us describe ${\cal G}_2$ explicitly as follows.

\vskip .1cm

\proclaim (4.5.6) Lemma. As a sheaf of sets,
$${\cal G}_2 \simeq \hbox{Hom}(T_M^{\otimes 2},T_M) \times
\hbox{Hom}(T_M^{\otimes
3},T_M).$$
Under this identification, the product of two sections $\varphi =
(\varphi_2, \varphi_3)$ and $\psi = (\psi_2, \psi_3)$, is given by
$$\eqalign{%
\varphi\psi &= (\varphi_2 + \psi_2, \varphi_3 + \psi_3 + \varphi_2
\circ (1 \otimes \psi_2 + \psi_2 \otimes 1)),\ {\rm and}\cr
\varphi^{-1} &= (-\varphi_2, -\varphi_3 + \varphi_2 \circ (1 \otimes
\varphi_2 + \varphi_2 \otimes 1))}.$$

\noindent {\sl Proof:} The identification is obtained first for the steaf of Lie
algebras corresponding to ${\cal G}_2$, by using the ${GL}(T_M)$-action and
then translated to ${\cal G}_2$ via the
exponential map.  To see that the multiplication has precisely this
form, it is enough to consider just one vector space $T$ and the group
$G_2$ of automorphisms of ${\cal A}ss^{\leq 3}(T^*)$ with identical
linear part.  In coordinates $x_1, \ldots, x_n$ an element of $G_2$ can
be written as
$$x_i \mapsto x_i + \sum_{j,k} a_i^{jk} x_j x_k + \sum_{p,q,r}
a_i^{pqr} x_p x_q x_r$$
where $\|a_i^{jk}\|$ form $\varphi_2 \in \hbox{Hom}(T^{\otimes 2}, T)$ and
$\|a_i^{pqr}\|$ form $\varphi_3 \in \hbox{Hom}(T^{\otimes 3},T)$.  Our
lemma is obtained by writing the composition of two such changes
of variables in the tensor notation.

\vskip .1cm

To prove (4.5.5), let $\{U_i\}$ be a Zariski open covering of ${\cal M}$ and
$\Phi = (\varphi_2^{(ij)})$ be a \v{C}ech 1-cocycle in
$\hbox{Hom}(T_M^{\otimes 2},T_M)$ representing $\alpha$.  To find the
2-cocycle $\psi = \delta(\Phi)$, we need to lift each $\varphi_2^{(ij)}$
to a section $\varphi^{(ij)}$ of ${\cal G}_2$ and set
$$\psi^{(ijk)} = (\varphi^{(ik)})^{-1} \varphi^{(jk)}
\varphi^{(ij)}.$$
Because of the above lemma, we can take $\varphi^{(ij)} =
(\varphi_2^{(ij)}, 0)$ and then an immediate calculation gives
$$\psi^{(ijk)} = (0, \varphi_2^{(jk)} \circ (1 \otimes \varphi^{(ij)}
+ \varphi^{(ij)} \otimes 1))$$
which is the Alexander-Whitney product of $\Phi$ and $(1 \otimes \Phi
+ \Phi \otimes 1)$.  Theorem is proved.

\vskip .1cm

The anti-associativity property becomes ordinary associativity under
suspension.

\proclaim (4.5.7) Corollary. Let $X$ be a 1-smooth thickening of ${\cal M}$
extendable to second order, and $\alpha = \alpha_{2,X}$.  Let ${\cal
A}$ be a sheaf of commutative ${\cal O}_M$-algebras.  Then the operation
$$H^i(M,{\cal A} \otimes T_M) \otimes H^j(M,{\cal A} \otimes T_M)
\mathop{\to}\limits^\cup H^{i+j}(M,{\cal A} \otimes T_M^{\otimes 2})
\mathop{\to}\limits^{(-1)^i\alpha} H^{i+j+1}(M,{\cal A} \otimes T_M)$$
makes $H^{\bullet -1}(M,{\cal A} \otimes T_M)$ into a graded associative
algebra.  When ${\cal A} = {\cal O}_M$, this algebra is graded commutative.

\noindent {\sl Proof.} The associativity follows by direct checking of signs,
using the Koszul rules.  The graded Lie algebra associated to
$H^{\bullet-1}(M,{\cal A} \otimes T_M)$ is clearly induced by the
symmetrization of $\alpha$, which is $\alpha_{2,M}$.  This graded Lie
algebra structure was studied in [Kap], where it was proved that for
${\cal A} = {\cal O}_M$ the bracket vanishes.

\vskip .3cm

\def\Hom{{\rm Hom}}

\noindent
{\bf (4.6) Operadic interpretation of the obstruction via Stasheff
polytopes.}
We now give a conceptual interpretation of the obstruction $\gamma_X$ from
Theorem 4.3.2, for arbitrary d. For this, we will use the language of
operads, see, e.g., 
 [GK1] [GK2],  and follow the conventions of [Kap], \S 3. In
particular, by $\Sigma{\cal P}$ we denote the suspension of a dg-operad
${\cal P}$.

Given a graded operad ${\cal P}$, an element $\xi\in {\cal P}(2)$ of
degree 0 will be called associative, if $\xi(\xi, 1) = \xi(1, \xi)$ in
${\cal P}(3)$.

Let ${ M}$ be a smooth algebraic variety, $T = T_M$. The sheaves
${\cal E}_{T}(n) = \Hom({T}^{\otimes n}, {T})$ form a
sheaf ${\cal E}_{T}$ of operads on ${\cal M}$. Thus the  ${\bb E}_{T}(n)
= {H}^{\bullet}({M}, {\cal E}_{T}(n)$ form a graded
operad ${\bb E}_{T}$. Given a 1-smooth thickening $X \supset M$, it
NC-Atiyah class $\alpha = \alpha_{2, {X}}$ is an element of ${\bb
E}_{T}(2)$ of degree 1. Thus the de-suspension
$\Sigma^{-1}(\alpha)\in  (\Sigma^{-1}{\bb E}_{T})(2)$ is an
element of degree 0. Theorem~4.5.3 can be formulated as follows.

\vskip .1cm

\proclaim (4.6.1) Corollary. $X$ can be extended to a
2-smooth
thickening if and only if $\Sigma^{-1}(\alpha)$ is an associative element.

We now construct a $dg$-model for ${\bb E}_{T}$. Fix an affine covering
$\{U_i\}$ of $M$, with intersections $U_{i_o\cdots i_m}$ affine as well. Let
$j_{i_o\cdots i_m} : U_{i_o\cdots i_m} \hookrightarrow {M}$ be the
embedding and
$$\check{\cal C}^m = \bigoplus_{i_o\cdots i_m} (j_{i_o\cdots i_m})_{\ast}
j_{i_o\cdots i_m}^{\ast} ({\cal O}_{M}), \quad \check{{C}}^m =
\Gamma({M}, \check{\cal C}^m),$$
so that $\check{{ C}}^m$ is the $\it m$-th term of the \v{C}ech
resolution of ${\cal O}_{M}$. The $\check{{\cal C}}^m$ form a
cosimplicial sheaf $\check{{\cal C}}^{\bullet}$ of commutative algebras
on $M$. Let ${\cal C}^{\bullet}$ be the sheaf of commutative $dg$-algebras
obtained from $\check{\cal C}^{\bullet}$ by the Thom-Sullivan
construction [HS]. For a quasicoherent sheaf ${\cal F}$ on $M$ let
 $${\cal C}^{\bullet}({\cal F}) = {\cal C}^{\bullet}\otimes {\cal F}, \quad
 C^{\bullet}({\cal F}) = \Gamma({M}, {\cal C}^{\bullet}({\cal F})).$$
The following is then standard.

\vskip .1cm

\proclaim  (4.6.2) Proposition. (a) ${\cal
C}^{\bullet}$ is flat over ${\cal O}_{M}$.\hfill\break
(b) For every ${\cal F}$ as above the ${\cal C}^m({\cal F})$ have
no higher cohomology. In particular, ${C}^{\bullet}({\cal F})$
calculates ${H}^{\bullet}({M}, {\cal F})$.\hfill\break
(c) There are natural morphisms of complexes
$${C}^{\bullet}({\cal F}_1)\otimes {C}^{\bullet}({\cal F}_2)
\rightarrow {C}^{\bullet}({\cal F}_1\otimes {\cal F}_2)$$
compatible with permutations of tensor factors.

It follows from (c) that the complexes ${
C}^{\bullet}({\cal E}_T({n}))$ form a $dg$-operad ${
C}^{\bullet}{\cal E}_{T}$ whose cohomology operad is ${\bb
E}_{T}$.

Each time when we have a $dg$-operad ${\cal P}$ and an associative element
$\xi\in  {H}^{0}{\cal P}(2)$, we can ask, following Stasheff
[Sta], about a hierarchy of finer properties on the level of cochains.
Namely, if ${m}_2\in  {\cal P}(2)^{0}$ is a cocycle representing
$\xi$, and ${m}_3\in  {\cal P}(3)^{-1}$ satisfies
$$\partial{m}_3 = {m}_2({m}_2, 1) - {m}_2(1, {m}_2)$$
(it exists because of the associativity of $\xi$), the element
$${m}_3(1, 1, {m}_2) + {m}_3({m}_2, 1, 1) -
{m}_2(1, {m}_3) - {m}_3(1, {m}_2, 1) -
{m}_2({m}_3, 1) \in  {\cal P}(4)^{-1}\leqno{(4.6.3)}$$
is a cocycle. If we visualize ${m}_2$ as a binary operation $({a,b})\mapsto
{ab}$ and ${m}_3$ as an edge (homotopy) connecting ${a(bc)}$ and
${(ab)c}$, then (4.6.3) represents the boundary of the MacLane-Stasheff
pentagon whose vertices correspond to the parenthesizings of ${abcd}$.
More generally, we have the following definition [Sta].

\vskip .1cm

\proclaim (4.6.4) Definition. Let ${\cal P}$ be a
$dg$-operad, $\xi\in {H}^{0}{\cal P}(2)$ an associative
element. We say that $\xi$ is an ${A}_{d}$-element if there
exist ${m}_i,\in {\cal P}(i)^{2-i}$, $2\le i \le {d}$
such that ${m}_2$ is a cocycle representing $\xi$, and for each
$\nu \in [2,{d}]$, we have
$$\partial {m}_\nu = \sum_{\mathop{%
r + {s} = \nu + 1}\limits_{\scriptstyle r, {s} \ge 2}}%
\enspace \sum_{1\ge k \ge r}(-1)^{{k}({s} + 1)} {m}_{r}
\mathop{\underbrace{(1,\dots, 1}}\limits_{{k} - 1},
{m}_{s}, 1,\ldots, 1).$$
Such a system $({m}_2,\ldots {m}_{d})$ will be called an ${
A}_{d}$-structure on $\xi$.

Geometrically, the expression on the right represents the boundary of the
Stasheff polytope, see [Sta].

\vskip .1cm

We apply this concept to ${\cal P} = \Sigma^{-1}{
C}^{\bullet}{\cal E}_{T}$.

\proclaim (4.6.5) Theorem. Let ${ M}$ be a smooth
algebraic
variety, ${X}\supset {M}$ a d-smooth thickening and
$\alpha\in {H}^1({M},~{\rm Hom}({T}_{M}^{\otimes 2},
{T}_{M}))$ the NC-Atiyah class of ${X}^{\le 1}$. Then:\hfill\break
{(a)} $\Sigma^{-1}(\alpha)$ is an ${A}_{d + 1}$-element.\hfill\break
{(b)} The image of the obstruction $\gamma_{X}$ in ${H}^2({M},~{\rm
Hom}({T}_{M}^{\otimes({d}+2)}, {T}_{M}))$, see
(4.4.5), is the cohomology class of the cocycle
$$\Sigma \left(\sum_{\mathop{%
r + {s} = d + 2}\limits_{\scriptstyle 1\le k \le r}}%
(-1)^{{ k}({s} + 1)} {m}_{r}
\mathop{\underbrace{(1,\dots, 1}}\limits_{{ k} - 1},
{m}_{s}, 1,\ldots, 1)\right)\in {C}^2%
\left({\rm Hom}({T}_{M}^{\otimes({d}+ 2)}, {T}_{M})\right)$$
where $({m}_2, \ldots, {m}_{{d}+1})$ is an appropriate ${
A}_{{d}+1}$-structure on $\Sigma ^{-1}(\alpha)$.

\noindent {\sl Proof:} It is an observation of Stasheff [Sta] that
an ${A}_{d}$-structure on a monoid $G$ is precisely the data needed to
build the classifying space of $G$ up to the $d$th level. We apply this to
the situation when $G$ is the shifted complex (${\cal C}^{\bullet -1}({\rm
T}_{M}),\partial$) of sheaves of ${\cal O}_{M}$-modules, the
operation is given by a 1-cocycle of ${m}_2$ representing $\alpha$ and
the skeleton of the classifying space is the sheaf of graded algebras ${\cal
C}^{\bullet}({\cal A}{ ss}^{\le d}(\Omega_{M}^1))$ with an
appropriate differential.

So we consider the sheaf ${\cal A}{ ss}^{\le ({d}+1)}(\Omega_{M}^1)$
of ${\cal O}_{M}$-algebras, filtered by the powers $\underline{\bf m}^i$.
Accordingly,
${\cal C}^{\bullet}({\cal A}{ ss}^{\le ({d}+1)} (\Omega_{M}^1))$
is a sheaf of graded, filtered ${\cal C}^{\bullet}({\cal O}_{M})$-algebras.

\proclaim (4.6.6) Proposition. ${A}_{{
d}+1}$-sequences $({m}_2,\ldots, {m}_{{d}+1})$, ${
m}_i\in  {C}^1~{\rm Hom}({T}_{M}^{\otimes(i)},~{T}_{
M})$, are in bijection with algebra differentials ${\cal D}$ of ${\cal
C}^{\bullet}({\cal A}{ ss}^{\le({d}+1)}(\Omega_{M}^1))$ of degree
1, satisfying the following conditions:\hfill\break
{(1)} ${\cal D}({ am}) = \partial ({a}){m} + (-1)^{\deg(a)}
{ a}\cdot {\cal D}({m}), \quad {\rm a}\in {\cal C}^{\bullet}({\cal
O}_{M})$, ${m}\in {\cal C}^{\bullet}({\cal A}{ ss}^{\le({\rm
d}+1)}(\Omega_{M}^1))$.\hfill\break
{(2)} ${\cal D}$ preserves the filtration and ${\cal D}\equiv
\partial\otimes 1$ on the associated graded factors.\hfill\break
{(3)} ${\cal D}^2 = 0$.

\noindent {\sl Proof:} An algebra derivation ${\cal D}$ satisfying
(1-2), is uniquely defined by its restriction on $\Omega_{M}^1$, which
gives morphisms $\tilde{m}_i : \Omega_{M}^1\rightarrow {\cal C}^1
((\Omega_{M}^1)^{\otimes i})$ or, equivalently, 
$${m}_i \in
{\cal C}^1({\rm Hom}({T}_{M}^{\otimes i},~{T}_{M})), i =
2,\ldots, {d} + 1.$$
  The condition that $({m}_2,\ldots,{m}_{{
d}+1})$ form an ${A}_{{d}+1}$-sequence, is equivalent to ${\cal D}^2
= 0$. Proposition is proved.

\vskip .1cm

\proclaim (4.6.7) Proposition. Given a derivation ${\cal D}$
of ${\cal C}^{\bullet}({\cal A}{ ss}^{\le ({d}+1)}(\Omega_{\rm
M}^1))$, the subsheaf ${\rm Ker} ({\cal D})$ $\i$  ${\cal C}^{0}({\cal A}{\rm
ss}^{\le ({d}+1)}(\Omega_{M}^1))$ is a sheaf of ${\cal O}_{\rm
M}$-algebras locally isomorphic to ${\cal O}_{M}\otimes({\bb C} {\<}
{ x}_1,\ldots, { x}_n {\>}/{\bf m}^{{d}+2})$ and equipped with an
identification $\underline {\bf m}/\underline {\bf m}^2\simeq \Omega_{M}^1$.

\noindent {\sl Proof:} Recalling the definition of ${\cal
C}^{\bullet}$ as the Thom-Sullivan construction applied to the \v{C}ech
complex, we find that ${\rm Ker}({\cal D})$ is obtained by gluing the
sheaves of algebras ${\cal A}{\rm ss}^{\le({d}+1)}(\Omega_{ U_i}^1)$
according to a system of descent data.

Conversely, given a sheaf ${\bb O}$ of ${\cal O}_{M}$-algebras as in
(4.6.7), we identify ${\bb O}|_{ U_i}
\tilde{{\displaystyle\rightarrow}}
{\cal A}{ ss}^{\le({d}+1)}$  $(\Omega_{ U_i}^1)$ by isomorphisms
identical on $\underline {\bf m}/\underline {\bf m}^2$ and get a derivation ${\cal D}$.
Proposition is proved.

\vskip .1cm

This proves part (a) of Theorem 4.6.5. Further, this proves that whenever
the cocycle in (b) is a coboundary, we can construct a sheaf ${\bb O}^{({\rm
d}+1)}$ on $M$ with ${\bb O}^{({d}+1)}/m^{{d}+2}\simeq { J}_{\rm
X}^{d}$. So the vanishing of $\varepsilon_{\ast}(\gamma_{X})$ is
equivalent to the vanishing of the cohomology class of the cocycle in (b).
Further, this observation can be upgraded to the equality of the two
cohomology classes, if we use once again the language of stacks. More
precisely, $\varepsilon_{\ast}(\gamma_{X})$ is represented, according to
(4.4.2)(d), by the stack of extensions
$$0\rightarrow(\Omega_{X}^1)^{\otimes({d}+2)}\rightarrow
{\bb O}^{({d}+1)}\rightarrow {\rm J}_{X}^{d}\rightarrow 0%
\leqno{(4.6.8)}$$
If we represent ${ J}_{X}^{d}$ as ${\rm Ker}({\cal D}) \subset
{\cal C}^{0}({\cal A}{ ss}^{\le ({d}+1)}(\Omega_{M}^1))$ with
${\cal D}$ given by $({m}_2,\ldots {m}_{{d}+1})$, then the
cocycle in (4.6.5)(b) is just $\tilde{{\cal D}}^2$, where $\tilde{{\cal D}}$
is the natural extension of ${\cal D}$ to ${C}^{\bullet}({\cal A}{
ss}^{\le ({d}+2)}(\Omega_{M}^1))$, obtained by taking $\tilde{{
m}}_{{d}+2} = 0$, see the proof of (4.6.6).  For any Zariski open ${
U}\subset {M}$ and any $\mu_{{d}+2}\in  {C}^1({\rm
Hom}({T}_{ U}^{\otimes({d}+2)},~{T}_{ U}))$ such that
$\partial\mu_{{d}+2}$ equals the restriction of the cocycle in
(4.6.5)(b) to $ U$, the sequence $({m}_2,\ldots, {m}_{{
d}+1}, \mu_{{d}+2})$ is an ${A}_{{d}+2}$-sequence and thus
determines an extension (4.6.8) over $ U$.

On the other hand, for any coherent sheaf ${\cal F}$ on ${M}$ and any
cocycle ${c}\in {C}^2 ({\cal F}) = \Gamma({M}, {\cal
C}^2({\cal F}))$ the class $[{c}]\in ~{H}^2({M}, {\cal F})$
can be represented by the following stack (${\cal F}$-gerbe) ${\cal S}_{\rm
c}$. The object of ${\cal S}_{c}({ U})$, ${ U}\subset {M}$, are
sections of $\partial^{-1}({c})|_{ U} \subset {\cal C}^1({\cal F})|_{\rm
U}$. A morphism from $\it b$ to $b^\prime$, where $b,
b^\prime\in ~{\rm Ob} ({\cal S}_{c}({ U}))$, are sections
$\alpha\in \Gamma({ U},~{\cal C}^{0}({\cal F}))$ such that
$\partial(a) = b - b^\prime$.

Now, taking for ${c}$ the cocycle from (4.6.5)(b), we find that
the above reasoning defines in fact an equivalence of ${\rm Hom}({\rm
T}_{M}^{\otimes({d}+2)},~{T}_{M})$-gerbes between
${\cal S}_{c}$ and the gerbe of extensions (4.6.8).
Theorem 4.6.5 is proved.

\vskip .2cm

\noindent {\bf (4.6.9) Remarks.} Let  $X$ be a $d$-smooth thickening of $M$.
 It is natural to call a (left) vector
bundle on $X$ a locally free sheaf $\cal E$ of left ${\cal O}_X$-modules.
Such a sheaf defines a vector bundle ${\cal E}_{ab}$ on $M$ in the usual
sense. Given a vector bundle
$E$ on $M$, the problem of extending $E$ to a bundle $\cal E$ on $X$ 
(with ${\cal E}_{ab}=E$) can be analyzed
in a way entirely parallel to the above. Thus, $E$ defines the
Atiyah class
$$\alpha_E\in H^1(M, {\rm Hom}(T\otimes E, E)),$$
see, e.g., [Kap]. The condition that $E$ extends to a bundle on $X^{\leq 1}$,
is: $$\alpha_E\circ (\alpha_X\otimes 1) = - \alpha_E\circ (1\otimes\alpha_E)
\quad {\rm in}\quad H^2(M, {\rm Hom}(T_M\otimes T_M\otimes E, E),$$
where $\alpha_X = \alpha_{2,X}$ is the NC-Atiyah class of $X^{\leq 1}$. After suspension,
this becomes the associativity condition for a (left) module over an algebra,
thus fitting nicely with (4.5.3) which gives the associativity condition
for an algebra. The obstructions to extending $E$ to bundles
on the further $X^{\leq i}$ can be interpreted in a similar way,
using the concept of a left $A_i$-module over an $A_i$-algebra. 
We leave this to the reader.

\vfill\eject

\centerline {\bf \S 5. Examples of NC-manifolds.}

\vskip 1cm

\noindent
{\bf (5.1) The projective space.} We define the
NC-manifold ${P}_{\rm left}^{n} = ({P}^{n}, {\cal
O}_{\rm left}^{\rm NC})$ (the left NC-projective space) by gluing ${
n}+1$ copies ${ U}_i, i=0,\ldots,{n}$, of ${A}_{\rm NC}^{
n}$. Namely, let ${x}_{ j}^{({ i})}$, $j\not= i$, be the
coordinates in ${ U}_i$, so that ${ U}_i = {\rm Spf}~{\bb C} \langle
{ x}_{ j}^{({ i})}\rangle_{\ldb ab\rdb}$. We now use the
noncommutative
variables $t_0,\ldots, {t}_{n}$ to relate the ${ x}_{ j}^{({
i})}$ together via
$${ x}_{ j}^{({ i})} = {t}_{ j} {t}_{ j}^{-1}%
\leqno{(5.1.1)}$$
which gives the transition functions
$${ x}_{ j}^{({ k})} = \left\{\eqalign{%
&{ x}_{ j}^{({ i})} ({ x}_{ k}^{({ i})})^{-1},\quad
  { j}\not= { i}\cr
&({ x}_{ k}^{({ i})})^{-1},\quad { j}= { i}\cr}\right.%
\leqno{(5.1.2)}$$
These expressions should be understood as commutator expansions, as
explained in Example (3.5.4)(a). One verifies immediately that the
functions (5.1.2) satisfy the cocycle conditions and thus define, by gluing,
an NC-manifold ${P}_{\rm left}^{n}$.

The right projective space ${P}_{\rm right}^{n} = ({P}^{n},
{\cal O}_{\rm right}^{\rm NC})$ is defined in a similar way, using the other
order of division. One verifies immediately that ${\cal O}_{\rm right}^{\rm
NC}$ is anti-isomorphic, as a sheaf of rings, to ${\cal O}_{\rm left}^{\rm
NC}$.

Further, we construct a sheaf ${\cal O}_{\rm left}^{\rm NC} (-1)$ of left
${\cal O}_{\rm left}^{\rm NC}$-modules, locally free and of rank 1, by using
the 1-cocycle of the covering $\{{ U}_{ i}\}$ in $({\cal O}_{\rm
left}^{NC})^{\ast}$ given by
$${\phi}_{{ i}_{ j}} = {t}_{ i} {t}_{ j}^{-1} = {
x}_{ j}^{({ i})} = ({x}_{ i}^{({ j})})^{-1}\in \Gamma
({ U}_{ i}\cap { U}_{ j}, ({\cal O}_{\rm left}^{\rm
NC})^{\ast}).$$
The dual sheaf ${\cal O}_{\rm left}^{\rm NC}(+1) = {\rm Hom}({\cal O}_{\rm
left}^{\rm NC}(-1), {\cal O}_{\rm left}^{\rm NC})$ is a sheaf of right ${\cal
O}_{\rm left}^{\rm NC}$-modules, locally free of rank 1. The sheaves ${\cal
O}_{\rm left}^{\rm NC}({m})$ for other ${m}$ do not make sense.

\vskip .2cm

\noindent
{\bf (5.1.3) Remark.} Thus our construction of
noncommutative projective spaces differs from the ``noncommutative
projectivization'' of graded algebras as developed in [AZ] [Ros] [VW]. In
our
setting, a sheaf of left ${\cal O}$-modules on ${P}_{\rm
left}^{n}$ does not give rise to any graded module since the ${\cal
O}({m})$, ${m}\not= \pm 1, 0$, are not defined. In fact, for a
graded NC-milpotent algebra ${ R}$ one can define a ringed space ${\rm Proj}(R)$
with the underlying topological space ${\rm Proj}({ R}_{ab})$ by
considering degree 0 elements in homogeneous localizations (cf. \S 1).
Applying this construction to the algebras ${\bb C} {\<} {t}_0,\ldots,{
t}_{n}{\>}/{ F}^{{d}+1}$ and passing to the inverse limit as ${
d}\rightarrow \infty$, we get a sheaf of algebras ${\cal O}^{\rm big}$ on
${P}^{n}$. This sheaf contains both ${\cal O}_{\rm left}^{\rm NC}$
and ${\cal O}_{\rm right}^{\rm NC}$ but is strictly bigger than any of them.
In particular, $({P}^{n}, {\cal O}^{\rm big})$ is not an
NC-manifold.

\vskip .3cm

\noindent
{\bf (5.2) Grassmannians.} Let $G(m,n)$ be the
Grassmannian of $m$-dimensional subspaces in ${\bb C}^{n}$. The
corresponding representable functor $h_0: {\cal C}om\rightarrow {\cal S}ets$ 
has the following explicit description. For a commutative
algebra $\Lambda$ the set $h_0(\Lambda)$ consists of all submodules $V
\subset \Lambda^n$ which are direct summands (hence projective) and
have rank m.

Now let us define a function $h:{\cal N}\rightarrow {\cal S}ets$ in a
similar way: $h(\Lambda)$ is the set of left submodules ${ V}\subset
\Lambda^{n}$ which are direct summands and whose rank (as
projective $\Lambda$-modules) is $m$.

\vskip .1cm

\proclaim  (5.2.1) Theorem. The functor $h$ is formally smooth and
commutes with arbitrary fiber products (2.3.3) in which at least one
of the $p_i$ is surjective.  Therefore (2.3.5) $h$ is represented by
an NC-smooth thickening $G(m,n)_{\rm left}^{\rm NC} \supset G(m,n)$.

\noindent {\sl Proof;} To see the formal smoothness, it is enough to prove that
$h(\Lambda') \to h(\Lambda)$ is surjective for any central extension
$$0 \to I \to \Lambda'\to \Lambda \to 0$$
of NC-nilpotent algebras.  Let us prove the following more general
fact.

\proclaim  (5.2.2) Lemma. In the above situation let ${P}'$ be
a projective (left) $\Lambda'$-module, ${P} = \Lambda \otimes
P'$ and $Q \subset {P}$ be a direct summand.  Then the set
of direct summands $Q' \subset P'$ such that $\Lambda
\mathop{\otimes}\limits_{\Lambda'} Q' = Q$, is a principal homogeneous
space over
$$I \mathop{\otimes}_{\Lambda_{ab}}
\Hom_{\Lambda_{ab}}(Q_{ab},P_{ab}/Q_{ab}).$$

\noindent {\sl Proof:} It is enough to prove this under extra assumptions which
can be satisfied by passing to localization.  Then we would have that
the direct summands $Q'$ form a sheaf on the Zariski topology of ${\rm
Spec}(\Lambda_{ab})$ which is a torsor over the coherent sheaf
corresponding to the $\Lambda_{ab}$-module specified above.  Since
$H^1$ of a coherent sheaf on an affine scheme vanishes, this torsor
has a global section and thus the set of all its global sections is a
principal homogeneous space as claimed.

Now, the extra assumptions we need, are:

\item{(1)}{$P'$ is free, $P' \simeq (\Lambda')^n$}

\item{(2)}{$Q$ is free, $Q \simeq \Lambda^m \subset \Lambda^n$.}

Let $e'_1, \dots, e'_n$ be the standard basis of $(\Lambda')^n$, and
$e_1, \dots, e_n$ the corresponding basis of $\Lambda^n$.  Since $Q =
\langle e_1, \dots, e_m\rangle$, an extension $Q'$ of $Q$ is spanned
by
$$e'_i + \sum_{j=m+1}^n b_{ij} e'_j, \qquad i=1,\dots,m, \quad b_{ij}
\in I.$$
So there is a bijection between all extensions of $Q$ and all matrices
$B = \|b_{ij}\| \in {\rm Mat}_{l,m-l}(I)$.  The set of such matrices
can be written invariantly as the tensor product 
\hfill\break
$I
\mathop{\otimes}\limits_{\Lambda_{ab}} \Hom_{\Lambda_{ab}}(Q_{ab},
P_{ab}/Q_{ab})$.  Even though the constructed identification of the
set of extensions with this tensor product is not canonical, the
corresponding structure of a principal homogeneous space is, as it can
be checked directly.  This proves Lemma 5.1.2 and the formal
smoothness of $h$.

\vskip .1cm

Let now $p_i: \Lambda_i \to \Lambda$, $i = 1,2$, be given and at least
one of the $p_i$ is surjective.  In addition to the natural map
$$j: h(\Lambda_1 \mathop{\times}\limits_{\Lambda} \Lambda_2) \to
h(\Lambda_1) \mathop{\times}\limits_{h(\Lambda)} h(\Lambda_2),$$
let us define a map
$$k: h(\Lambda_1) \mathop{\times}\limits_{h(\Lambda)} h(\Lambda_2) \to
h(\Lambda_1 \mathop{\times}\limits_{\Lambda} \Lambda_2)$$
to associate to a pair $({ V}_1 \subset \Lambda_1^n, { V}_2 \subset
\Lambda_2^n)$ with $\Lambda \mathop{\times}\limits_{\Lambda_1} V_1 =
\Lambda \mathop{\times}\limits_{\Lambda_2} V_2 = V \subset \Lambda^n$,
the submodule $\tilde{V} = V_1 \mathop{\times}\limits_V V_2$ in
$\Lambda_1 \mathop{\times}_{\Lambda} \Lambda_2$.  We need only to
prove the lemma below.  Then one verifies at once that $j$ and $k$ are
mutually inverse.

\proclaim (5.2.3) Lemma.  $\tilde{V}$ is a direct summand of
$\tilde{\Lambda} = \Lambda_1 \mathop{\times}\limits_{\Lambda}
\Lambda_2$.

\noindent {\sl Proof:} If $\Lambda$ is any ring and $V \subset \Lambda^n$ is a
left submodule, then $\Lambda$ is a direct summand iff $\Lambda^n/V$
is projective.  So the set of direct summands in $\Lambda^n$
identified with the set of surjections $s: \Lambda^n \to P$ with $P$
projective, modulo isomorphisms of projective modules.

Applying this to our particular situation, we have surjections $s_i:
\Lambda_i^n \to P_i$, $s: \Lambda^n \to P$ and isomorphisms $\Lambda
\mathop{\otimes}\limits_{\Lambda_i} P_i {\to} P$
identifying $\Lambda \otimes s_i$ with $s$.  In this case $\tilde{P} =
P_1 \mathop{\times}\limits_P P_2$ is a projective
$\tilde{\Lambda}$-module ([Mil], Th. 2.1) and we have a surjection
$\tilde{s}: \tilde{\Lambda}^n \to \tilde{P}$ with kernel $\tilde{V}$.
This proves the lemma and Theorem 5.2.1.

\vskip.2cm

\noindent {\bf (5.2.4) Remarks.}  One can also construct $G(m,n)_{\rm
left}^{\rm NC}$ explicitly, by gluing ${n \choose m}$ copies of
$A_{\rm NC}^{m(n-m)}$.  For this, all we need is to mimic the standard
formulas [GH] and interpret the entries of the inverse
matrices involved, as commutator series, cf. (3.5.4)(c).  This
approach is closely related to the ``noncommutative Grassmannian''
construction of Gelfand and Retakh [GR].  In fact, $G(m,n)_{\rm
left}^{\rm NC}$ can be seen as the formal completion of their
construction along the commutative points.

\vskip .3cm

\noindent {\bf (5.3) Flag varieties.} Fix $1 \leq i_1-1 < \dots < i_k \leq n$
and let $F(i_1, \dots, i_k, {\bb C}^n)$ be the variety of flags ${ V}_1
\subset \dots \subset V_k \subset {\bb C}^n$, $\dim (V_{\nu}) =
i_{\nu}$.

\proclaim (5.3.1) Theorem. $F(i_1, \dots, i_k, {\bb C}^n)$ has a
natural NC-smooth thickening $F(i_1, \dots, i_k, {\bb C}^n)_{\rm
left}^{\rm NC}$, and the action of $GL_n({\bb C})$ extends to this
thickening.

\noindent {\sl Proof:} As in (5.2), we extend the functor represented by $F(i_1,
\dots, i_k, {\bb C}^n)$, to a functor $h: {\cal N} \to Sets$ with
$h(\Lambda)$ being the set of all flags of left submodules $V_1
\subset \dots \subset V_k \subset \Lambda^n$ which are direct summands
and such that the rank of $V_{\nu}$ is $i_{\nu}$.  The fact that $h$
is formally smooth and commutes with fiber products such as in (5.2.1)
is verified in a similar way, using Lemma 5.2.2 (plus induction) and
Lemma 5.2.3.

\vskip .3cm

\noindent {\bf (5.4) Moduli of vector bundles.}  Let ${ Z}$ be a
projective algebraic variety.  For any variety $B$ let $\rho: B \times
{ Z} \to B$ denote the projection.

\proclaim (5.4.1) Definition. A versal family of vector bundles on
$ Z$ is a pair $(B, E)$ where $B$ is a smooth algebraic variety,
$E$ a vector bundle on $B \times  Z$ with the following
properties:\hfill\break
{(a)}{$R^0\rho_*\ {\rm End}(E) = {\cal O}_B$.}\hfill\break
{(b)}{The Kodaira-Spencer map $\kappa: T_B \to R^1\rho_*\
{\rm End}(E)$ is an isomorphism.}\hfill\break
{(c)}{$R^2\rho_*\ {\rm End}(E) = 0$.}

Many (etale) open charts of moduli spaces of stable vector bundles are
bases of versal families.

\proclaim  (5.4.2) Theorem. If $(B, E)$ is a versal family of vector
bundles, then $B$ admits a canonical NC-smooth thickening.

\noindent {\sl Proof:} We construct a functor $h: {\cal N} \to {\cal }Sets$.  Let
$\Lambda$ be a NC-nilpotent algebra, $\Lambda_{ab}^0$ be the quotient
of $\Lambda_{ab}$ by the ideal of nilpotent elements, $X_{ab}^0 = {\rm
Spec}(\Lambda_{ab}^0)$.  The schemes $X_{ab}^0$ and $X_{ab} = {\rm
Spec}(\Lambda)$ having the same underlying topological space, we can
regard $X_{ab}^0$ as the underlying space for ${\rm Spec}(\Lambda)$,
i.e., write $X = {\rm Spec}(\Lambda) = (X_{ab}^0, {\cal O}_X)$.  Let $p_1,
p_2$ be the projection of $X_{ab}^0 \times { Z}$ to the factors.
Denote by $X \times { Z}$ the NC-scheme with underlying space
$X_{ab}^0 \times { Z}$ and the sheaf of algebras ${\cal O}_{X \times
{ Z}} = p_1^* {\cal O}_X \otimes p_2^* {\cal O}_{{ Z}}$.  Its
abelianization
is ${\cal O}_{X_{ab} \times { Z}}$.

Let ${\cal C}_{\Lambda}$ be the category whose objects are the
following sets of data:

\item{(1)}{A morphism of schemes $f: X_{ab}^0 \to B$}

\item{(2)}{A locally free sheaf $\cal E$ of left ${\cal O}_{X \times
{
Z}}$-modules.}

\item{(3)}{An isomorphism $\varphi: {\cal O}_{X_{ab}^0 \times { Z}}
\otimes {\cal E} \mathop{\to}\limits^{\sim} (f \times Id)^* E$.}

\noindent
A morphism $(f_1, \varepsilon_1, \varphi_1) \to (f_2, \varepsilon_2,
\varphi_2)$ in ${\cal C}_{\Lambda}$ exists only if $f_1 = f_2$ and
consists of an isomorphism $\varepsilon_1 \to \varepsilon_2$ commuting with
the $\varphi_2$.

\vskip .1cm

We define $h(\Lambda)$ to be the set of isomorphism classes of objects
of ${\cal C}_{\Lambda}$.  Our statement now will follow from (2.3.5)
and the next proposition.

\vskip .1cm

\proclaim (5.4.3) Proposition. (a) The restriction of $h$ to Com
coincides with $h_B$, the functor represented by $B$. \hfill\break
(b) $h$ is formally smooth.\hfill\break
(c) $h$ commutes with any fiber products (2.3.3) in which one of the
$p_i$ is surjective.

\noindent {\sl Proof:}  We first concentrate on (b).  It is enough to prove the
following.

\proclaim (5.4.4) Lemma. If
$$0 \to I \to \Lambda' \mathop{\to}\limits^p \Lambda \to 0$$
is a central extension of NC-nilpotent algebras and $[f, \varepsilon,
\varphi] \in h(\Lambda)$, then $h(p)^{-1}[f, \varepsilon, \varphi]$ is a
principal homogeneous space over $I
\mathop{\otimes}\limits_{\Lambda_{ab}} f^* R^1\rho_*\ {\rm End}(E)$.

\noindent {\sl Proof:} As in the proof of (5.2.2), it is enough to work locally
on $X_{ab}$, i.e., assume that $\Lambda_{ab}$ and $\Lambda$ are local
rings.  Then so is $\Lambda'$.  The set of isomorphism classes of
locally free sheaves of left ${\cal O}_{X \times Z}$-modules of rank $r$ is
in this case identified with $H^1(Z, GL_r({\cal O}_Z \otimes \Lambda))$, and
similarly for $\Lambda'$.  We have an exact sequence
$$1 \to {\cal O}_Z \times {\rm Mat}_r(I) \to GL_r(O_Z \otimes \Lambda')
\mathop{\to}\limits^{\tilde{p}} GL_r({\cal O}_Z \otimes \Lambda) \to 1$$
of sheaves of groups on $Z$.  An element $[f, {\cal E}, \varphi] \in
h(\Lambda)$ is represented by some $c \in H^1(Z, GL_r(O_Z \otimes
\Lambda))$, while
$$h(p)^{-1}[f, {\cal E}, \varphi] \cong \tilde{p}^{-1}(c) \subset
H^1(Z, GL_r({\cal O}_Z \otimes \Lambda')).$$
Our statement now follows from the general formalism of long exact
sequences of non-Abelian cohomology, as recalled, e.g. in [Man],
\S(2.6.8-9).  Lemma is proved.

\vskip .1cm

Thus we proved part (b) of (5.4.3).  To prove part (a), note that by
construction, $h = h_B$ on the category of commutative algebras
without nilpotents.  In addition, we have a natural transformation $u:
h_B \to h|_{\rm Com}$, which sends a morphism $\tilde{f}: {\rm
Spec}(\Lambda) \to B$ (with $\Lambda$ commutative) to the sheaf
$(\tilde{f} \times Id)^* E$.  We claim that $u$ is an isomorphism.  To
show this, it is enough to proceed by induction in the degree of
nilpotency of the ideal ${\rm Ker}(\Lambda \to \Lambda_0)$, and consider
a square-zero extension of commutative algebras, such as in (5.4.4) for
which $u_{\Lambda}: \Hom({\rm Spec}(\Lambda), B) \to h(\Lambda)$ is an
isomorphism.  Given $\tilde{f} \in h_B(\Lambda)$, the set
$h_B(p)^{-1}(\tilde{f})$ is a principal homogeneous space over $I
\mathop{\otimes}\limits_{\Lambda} \tilde{f}^* T_B$, while by (5.4.4),
$h(p)^{-1}(\tilde{f})$ is a principal homogeneous space over $I
\mathop{\otimes}\limits_{\Lambda} R^1 \rho_*({\rm End}(E))$.  Since
the Kodaira-Spencer map $\kappa: T_B \to R_1\rho_*({\rm End}(E))$
is an isomorphism by (5.4.1)(b), and $u_{\Lambda'}$ is compatible with
$\kappa$, it is a bijection.

\vskip .1cm

Finally, we prove part (c) of (5.4.3).  Given any diagram of algebras
as in (2.3.3), we define the map
$$k: h(\Lambda_1) \mathop{\times}\limits_{h(\Lambda)} h(\Lambda_2) \to
h(\Lambda_{12})$$
as follows.  Given $f_i: {\rm Spec}(\Lambda_{i,ab}^0) \to B$
coinciding on ${\rm Spec}(\Lambda_{12,ab}^0)$ and locally free sheaves
$\varepsilon_i$ of ${\cal O}_{X_i \times Z}$-modules, isomorphic upon the
pullback to ${\cal O}_{X_{12}} \times Z$, we define
$$k([f_1, {\cal E}_1, \varphi_1], [f_2, {\cal E}_2, \varphi_2]) =
[f_{12}, {\cal E}_{12}, \varphi_{12}]$$
where $f_{12} = f_1 \mathop{\times}\limits_f f_2$, ${\cal E}_{12} =
{\cal E}_1 \mathop{\times}\limits_{\cal E} {\cal E}_2$, and $f:
{\rm
Spec}(\Lambda_{12,ab}^0) \to B$, ${\cal E} = {\cal O}_{X_{12} \times Z}
\mathop{\otimes}\limits_{{\cal O}_{X_i} \times Z} {\cal E}_i$ are the
common
restrictions of $f_i$, ${\cal E}_i$ on $\Lambda$.  The sheaf ${\cal
E}_{12}$ is locally free by Theorem 2.1 from [Mil] already used in the
proof of (5.2.3).  Having constructed $k$, it is clear that it is the
inverse of the map $j$ from (2.3.4).  So both are bijections.

This completes the proof of Proposition 5.4.3 and Theorem 5.4.2.

\vskip .2cm

\noindent {\bf (5.4.6) Remark.}  The reason for all the examples considered
in this section, to possess natural NC-thickenings, can be traced to
the fact that he group $GL_n(\Lambda)$ can be defined for any
associative and not necessarily commutative algebra $\Lambda$.  This
can be expressed by saying that the natural NC-thickening $GL_n({\bb
C})^{\rm NC} \supset GL_n({\bb C})$ induced by the embedding
$GL_n({\bb C}) \subset A^{n^2} \subset A_{\rm NC}^{n^2}$, is a group
object in the category of NC-manifolds.

\vfill\eject

\centerline {\bf References}

\vskip 1cm

\noindent {[AZ]} M. Artin, J. J. Zhang, Noncommutative projective schemes, {\it
Adv. in Math.} 109(1994), 228-287.

\vskip .1cm

\noindent {[BE]} N. Blackburn, L. Evens, Schur multipliers of $p$-groups, {\it J.
Reine und Angew. Math.} 309(1979), 100-113.

\vskip .1cm

\noindent {[Bry]} J.-L. Brylinski, Loop spaces, characteristic classes and
geometric quantization, Birkh\"auser, Boston, 1993.

\vskip .1cm

\noindent {[Co]} A. Connes, Noncommutative geometry, Academic Press, 1994.

\vskip .1cm

\noindent {[CQ]} J. Cuntz, D. Quillen, Algebra extensions and nonsingularity,
{\it J. AMS,} 8(1995), 251-289.
\vskip .1cm

\noindent 
{[Eis]} D. Eisenbud, Commutative algebra with a view towards algebraic
geometry, Springer-Verlag, 1995.

\vskip .1cm

\noindent {[Ev]} L. Evens, Terminal $p$-groups, {\it Illinois J. Math.} 12(1968),
682-699.

\vskip .1cm

\noindent {[Fe]} R. Feynman, An operator calculus having applications in  quantum
electrodynamics, {\it Phys. Rev.} (2), 84(1951), 108-128.

\vskip .1cm

\noindent [GH] P.Griffiths, J. Harris, Principles of algebraic geometry,
J. Wiley, 1978. 

\vskip .1cm

\noindent {[GK1]} E. Getzler, M. Kapranov, Cyclic operads  and cyclic homology,
in ``Geometry, topology and physics for R. Bott'' (S. T. Yau, Ed.), P.
167-201, International Press, Cambridge MA, 1995.

\vskip .1cm

\noindent {[GK2]} E. Getzler, M. Kapranov, Modular operads, {\it Compositio Math.}
110(1998), 65-126.

\vskip .1cm

\noindent {[GR]} I. M. Gelfand, V. S. Retakh, Quasideterminants I, preprint
alg-geom/9705026, to appear in  {\it Selecta Math.}

\vskip .1cm

\noindent {[HH]} J. W. Helton, R. E. Howe, Traces of commutators of integral
operators, {\it Acta Math.} 135(1975), 271-305.

\vskip .1cm

\noindent {[HS]} V. S. Hinich, V. V. Schechtman, On homotopy limit of homotopy
algebras, in ``$K$-theory, arithmetic and geometry'' (Y. I. Manin Ed.), p.
240-264, Lecture Notes in Math. 1289, Springer-Verlag, 1987.

\vskip .1cm

\noindent {[Kap]} M. Kapranov, Rozansky-Witten invariants in Atiyah classes,
preprint alg-geom \# 9704009, to appear in {\it Compositio Math.}

\vskip .1cm

\noindent {[KM]} M. V. Karasev, V. P. Maslov, Nonlinear Poisson brackets,
geometry and quantization, {\it Amer. Math. Soc.} 1993.

\vskip .1cm

\noindent [Ko] M. Kontsevich, Formal non-commutative symplectic geometry,
in: ``Gelfand Mathematical Seminars 1990-92" (L.Corwin, I. Gelfand, J. Lepowsky Eds.),
173-187, Birkhauser, Boston, 1993.

\vskip .1cm

\noindent {[L]} J.-L. Loday, Cyclic homology, Springer-Verlag, 1992.

\vskip .1cm

\noindent {[MKS]} W. Magnus, A. Karrass, D. Solitar, Combinatorial group theory,
Dover Publ. 1976.

\vskip .1cm

\noindent {[Mac]} I. Macdonald, Symmetric functions and Hall polynomials, Oxford,
Clarendon Press, 1995.

\vskip .1cm

\noindent {[Man]} Y. I. Manin, Gauge field theory and complex geometry,
Springer-Verlag, 1991.

\vskip .1cm

\noindent {[Mil]} J. Milnor, Introduction to algebraic $K$-theory, Princeton
Univ. Press, 1971.

\vskip .1cm

\noindent {[Mas]} V. P. Maslov, Operational methods, Mir Publ. Moscow 1968.

\vskip .1cm

\noindent {[Ros]} A. L. Rosenberg, Noncommutative algebraic geometry and
representations of quantum groups, Kluwer Publ. 1995.

\vskip .1cm

\noindent {[Sche]} W. Schelter, Smooth algebras, {\it J. of Algebra,} 103(1986),
677-685.

\vskip .1cm

\noindent {[Schl]} M. Schlessinger, Functors on Artin  rings, {\it Trans. AMS,}
130(1968), 208-222.

\vskip .1cm

\noindent {[Sig]} S. Sigg, Laplacian and  homology of free two-step nilpotent Lie
algebras, {\it J. Algebra,} 185(1996), 144-161.

\vskip .1cm

\noindent {[Sta]} J. D. Stasheff, Homotopy associativity of $H$-spaces I, II,
{\it Trans. AMS,} 108(1963), 275-292.

\vskip .1cm

\noindent {[Ste]} B. Steinstrom, Rings of quotients, Springer-Verlag 1975.

\vskip .1cm

\noindent {[VV]} F. Van Oystaeyen, A. H. Verschoren, Non-commutative algebraic
geometry, Lecture Notes in Math. 887, Springer-Verlag 1981.

\vskip .1cm

\noindent {[VW]} F. Van Oystaeyen, L. Willaert, Cohomology of schematic algebras,
{\it J. Algebra,} 185(1996), 74-84.

\vskip 2cm

\noindent 
{\sl Author's address: Department of mathematics, Northwestern University,
Evanston IL 60208 USA, email: \hfill\break
kapranov@math.nwu.edu

\end